\documentclass[11pt]{amsart}
\usepackage[utf8]{inputenc}

\usepackage{amssymb, mathtools}
\usepackage[margin=1in]{geometry}
\usepackage{hyperref}
\usepackage{enumerate}
\usepackage{tikz-cd}

\newcommand{\bigo}{\mathcal{O}}
\newcommand{\R}{\mathbb{R}}
\newcommand{\Z}{\mathbb{Z}}
\newcommand{\N}{\mathbb{N}}
\newcommand{\Q}{\mathbb{Q}}
\newcommand{\C}{\mathbb{C}}
\newcommand{\A}{\mathbb{A}}
\newcommand{\p}{\mathfrak{p}}
\newcommand{\q}{\mathfrak{q}}
\renewcommand{\a}{\mathfrak{a}}
\newcommand{\trd}[1]{\mathrm{Trd}\left(#1\right)}
\newcommand{\nrd}[1]{\mathrm{Nrd}\left(#1\right)}
\newcommand{\GL}{\mathrm{GL}}

\newtheorem{theorem}{Theorem}[section]
\newtheorem{proposition}[theorem]{Proposition}
\newtheorem{corollary}[theorem]{Corollary}
\newtheorem{lemma}[theorem]{Lemma}

\theoremstyle{definition}
\newtheorem{definition}[theorem]{Definition}

\theoremstyle{remark}
\newtheorem{remark}[theorem]{Remark}
\newtheorem{remarks}[theorem]{Remarks}

\address{Department of Mathematics, University of Virginia, Charlottesville, VA 22904, USA}
\email{zha5sb@virginia.edu}
\urladdr{https://sites.google.com/view/maxsg/}


\allowdisplaybreaks

\begin{document}

\title{The geodesic restriction problem for arithmetic spherical harmonics}
\author{Maximiliano Sanchez Garza}
\date{}

\begin{abstract}
Given a Riemannian manifold $M$ and an $L^2$-normalized Laplacian eigenfunction $\psi$ on $M$ with eigenvalue $\lambda^2$, a general problem in analysis is to understand how the mass of $\psi$ distributes around $M$. There are different ways to attack this problem. One of them is to analyze the $L^p$-norm of $\psi$ restricted to a submanifold of $M$. Here, we concentrate on the case $M=S^2$, $p=2$, and we restrict to geodesics of the sphere. Burq, G\'erard, and Tzvetkov showed, for $\gamma$ a geodesic of $S^2$ (and indeed for more general surfaces), that $||\psi|_{\gamma}||_{L^2} \ll \lambda^{1/4}$ and that this bound is optimal in general. In this paper, we specialize to the case in which $\psi$ is an eigenfunction of all the Hecke operators on the sphere and consider the set of geodesics $\mathcal{C}_{D}$ of $S^2$ associated to fundamental discriminants $D<0$. By combining approaches of Ali and Magee, we improve the previous upper bound to $||\psi|_{\mathcal{C}_{D}}||_{L^2} \ll_{D,\varepsilon} \lambda^{\varepsilon}$ for any $\varepsilon>0$, which is essentially sharp.
\end{abstract}

\maketitle

\section{Introduction}\label{Introduction}

Let $(M,g)$ be a Riemannian manifold and let $\Delta$ be the Laplace operator associated to $g$. Let $\psi$ be a Laplacian eigenfunction on $M$ satisfying $\Delta \psi = \lambda^2 \psi$. An interesting problem in analysis is to study how much of the mass of $\psi$ concentrates in small subsets of $M$. There are several ways to approach this question; see the introduction of \cite{BGTGeneralBound} for some of these approaches and references. One of them is to take $\psi$ to be $L^2$-normalized and bound $||\psi|_{M'}||_{L^p}$, where $M'$ is a submanifold of $M$ and $2\leq p\leq \infty$. 

This approach was first studied by Reznikov in \cite{ReznikovInitialStudy} in the particular case of $M$ being a compact hyperbolic Riemannian surface, $p=2$, and restricting $\psi$ to a closed geodesic $\gamma$, where he showed that 
\begin{equation}
    ||\psi|_{\gamma}||_{L^2} \ll \lambda^{\frac{1}{4}}. \label{BGT-bound}
\end{equation}
Later, in \cite[Thm. 3]{BGTGeneralBound}, Burq, G\'erard, and Tzvetkov generalized Reznikov's results to compact Riemannian manifolds, where the bounds are now dependent on $p$ and the dimensions of $M$ and $M'$, and showed that these are sharp in general. Sarnak, in his letter to Reznikov in \cite[p. 1]{LetterSarnakReznikov}, calls bound \eqref{BGT-bound} and the bounds obtained by Burq, G\'erard, and Tzvetkov the \textit{convexity bounds}\footnote{The reason for this name is that, in the case of $M'=M$ being a compact arithmetic hyperbolic surface and $p=\infty$, there is a strong connection between our problem and deep conjectures in number theory such as the Ramanujan conjecture for Maass forms \cite[\S 5.1]{IwaniecKowalski} and the generalized Lindel{\"o}f hypothesis \cite[\S 5.2]{IwaniecKowalski}. This connection is described in \cite[\S 4]{SarnakAQC}.}.

Even though the bound \eqref{BGT-bound} is sharp in general, it is expected to not be optimal in the case where $M$ has constant curvature $-1$ since the geodesic flow is ergodic and quantum ergodicity holds (see \cite{Shnirelman}, \cite{Verdiere}, \cite{Zelditch}). The only power-saving progress in this setting is when $M$ is an arithmetic quotient of the complex upper half-plane, such as $M=\mathrm{SL}_{2}(\mathbb{Z}) \backslash \mathbb{H}$,\footnote{See \cite[\S 14.1]{IwaniecKowalski} for the definition of the action of $\mathrm{SL}_{2}(\mathbb{Z})$ on $\mathbb{H}$.} and $\psi$ is arithmetic, i.e.~a simultaneous eigenfunction of $\Delta$ and number-theoretic operators called the Hecke operators\footnote{See \cite[\S 14.6]{IwaniecKowalski}.}. The first improvement to the convexity bound in this setting is due to Iwaniec and Sarnak, where they showed in \cite[Thm. 0.1]{IwaniecSarnak} that $$ ||\psi||_{L^{\infty}} \ll_{\varepsilon} \lambda^{\frac{5}{12}+\varepsilon} $$ for any $\varepsilon>0$, improving the convexity bound $||\psi||_{L^{\infty}} \ll \lambda^{\frac{1}{2}}$. More recently, in \cite{DanaMScThesis}, Ali improved \eqref{BGT-bound} by showing that for any $\varepsilon>0$ we have 
\begin{equation}
    ||\psi|_{\mathcal{C}_D}||_{L^2} \ll_{D, \varepsilon} \lambda^{\theta + \varepsilon} \label{AbouAli-bound}
\end{equation}
where $D>0$ is a fundamental discriminant,\footnote{An integer $D\neq 1$ is a \textit{fundamental discriminant} if $D$ is not divisible by the square of any odd prime and it satisfies either $D\equiv 1 \pmod{4}$ or $D\equiv 8, 12 \pmod{16}$.} $\mathcal{C}_D$ is the union of all closed geodesics in $M=\mathrm{SL}_2(\Z) \backslash \mathbb{H}$ associated to $D$,\footnote{These closed geodesics are constructed from binary quadratic forms $ax^2 + bxy + cy^2$ of discriminant $ b^2 - 4ac = D$. An important feature of these closed geodesics, as shown by Duke in \cite{Duke}, is that the set $\mathcal{C}_{D}$ equidistributes on the surface $\mathrm{SL}_{2}(\mathbb{Z}) \backslash \mathbb{H}$ as $D \rightarrow \infty$. For more on this construction and Duke's theorem, see \cite{HarcosDukeThm}.} and $\theta \geq 0$ is any bound towards the Ramanujan Conjecture for Maass forms. One can take $\theta = \frac{7}{64}$ by results of Kim and Sarnak in \cite[Appendix 2, Prop. 2]{KimSarnak}. It is worth noting that, assuming the Ramanujan conjecture, we could take $\theta=0$. 

\begin{figure}[!htb]
    \centering
    \includegraphics[width=0.6\linewidth]{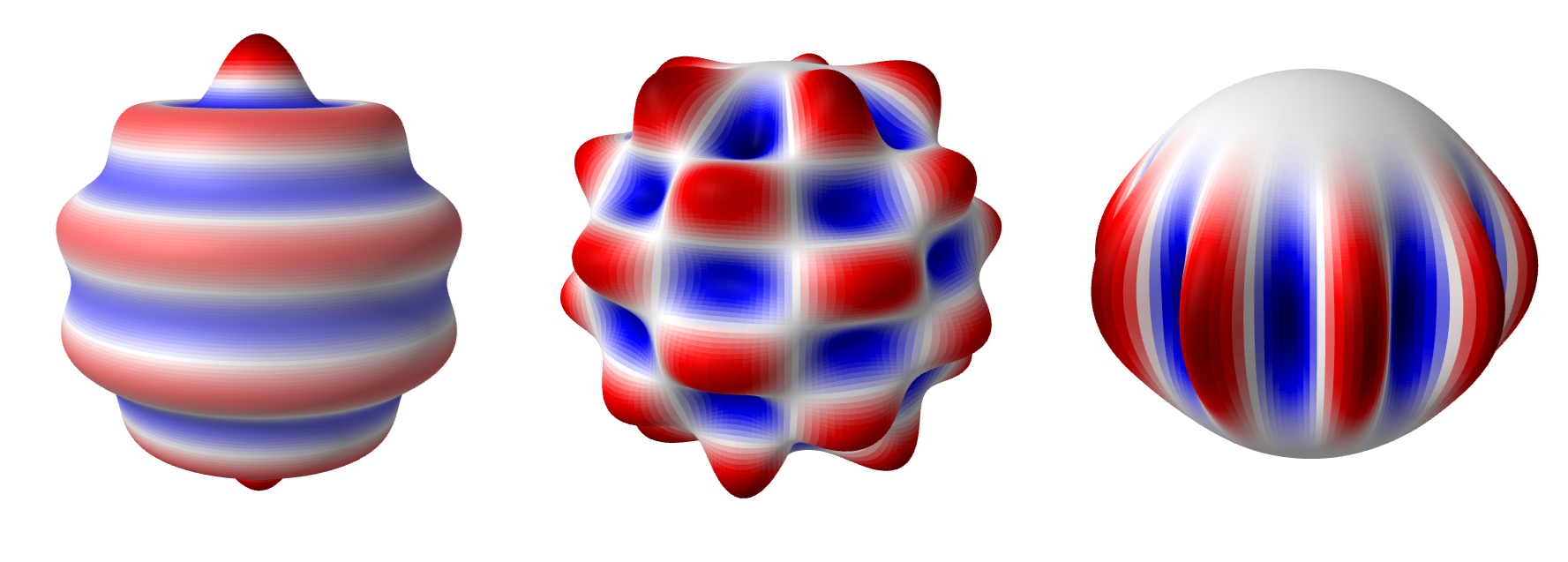}
    \caption{From left to right, graphic illustrations of zonal, tesseral, and sectorial spherical harmonics. If $f(\xi)$ is the spherical harmonic in consideration, the red zones represent positive values of $\mathrm{Re}(f(\xi))$, whereas the blue zones represent negative values of $\mathrm{Re}(f(\xi))$. The intensity of the color represents $|f(\xi)|$.}
    \label{fig1}
\end{figure}

We concentrate on the case of the two-dimensional sphere $S^2$ (or rather, a particular quotient of it, as we shall describe below), which has constant curvature $1$. We give it the round metric and associated measure, which gives it the structure of a Riemannian manifold. Let $\Delta$ be its Laplace operator. We know that the Hilbert space $L^2(S^2)$ has a complete spectral resolution with respect to $\Delta$. This spectral resolution is completely explicit in terms of the eigenfunctions of $\Delta$, called spherical harmonics; see \cite[Thm. 2.9]{SpectrumLaplacian}, \cite[p. 140]{SteinWeiss}, \cite[Th. 2.34, Prop. 3.5]{SpheHar}. More explicitly, if $$\{0 = \lambda_0 < \lambda_1 < \lambda_2 < \cdots \}$$ is the spectrum of $\Delta$, then $\lambda_\ell = \ell(\ell+1)$ and its eigenspace, denoted by $H_\ell$, is the space of harmonic homogeneous polynomials in three variables of degree $\ell$ restricted to the sphere $S^2$, which are precisely the spherical harmonics. 

Let $\psi \in H_\ell$ and let $\gamma$ be a geodesic of $S^2$, which is to say that $\gamma$ is a ``great circle'' of $S^2$, i.e.~the equator of $S^2$ when we choose a point of $S^2$ to be the north pole. We know $||\psi|_{\gamma}||_{L^2}$ satisfies the bound \eqref{BGT-bound}. In fact, for arbitrary $\psi$, the bound \eqref{BGT-bound} is sharp; see \cite[Thm. 1]{BGTGeneralBound}. The reason for this is that, at least for certain $\psi$, there are high concentrations of mass in specific regions of the sphere (see Figure \ref{fig1}). 

Although $S^2$ has constant curvature $1$, we expect improvements on the bound \eqref{BGT-bound} for specific $\psi$. This is due to the fact that a basis of Laplacian eigenfunctions for $S^2$ can also be diagonalized further with respect to Hecke operators, thus making the basis elements arithmetic. These operators were originally defined by Eichler in \cite{EichlerHeckeOp} using Brandt matrices, but they can be described explicitly using the action of the Hamiltonian quaternions on $S^2$ via rotations (see \cite[\S 3.1]{ConwaySmith}). 

More precisely, let $$B(\mathbb{R}) = \{ a + bi + cj + dk \: | \: a,b,c,d\in \mathbb{R}, \: i^2=j^2=k^2=-1, \: ij=-ji=k \}$$ be the Hamiltonian quaternions and $B^{\times}(\mathbb{R})$ be its group of units. For $\alpha = a+bi+cj+dk \in B(\mathbb{R})$, we let $$\mathrm{Trd}(\alpha) := 2a \in \mathbb{R} \: \text{ and } \: \mathrm{Nrd}(\alpha) := a^2+b^2+c^2+d^2 \in \mathbb{R}$$ denote its trace and norm, respectively. We identify $S^2$ with the elements of $B(\mathbb{R})$ of trace $0$ and norm $1$, and hence we have an action of $B^{\times}(\mathbb{R})\ni \alpha$ on $S^2 \ni z$ by conjugation, i.e.~$\alpha \cdot z := \alpha z \alpha^{-1}$. Consequently, there is an induced action on $L^2(S^2)$, which is the Koopman representation of $B^{\times}(\R)$ in $L^2(S^2)$. 

Let $\bigo\subset B(\R)$ denote the maximal order of Hurwitz integers, i.e. $$\bigo = \left\{ \left. \frac{a+bi+cj+dk}{2} \: \right| \: a,b,c,d\in \mathbb{Z}, \: a\equiv b\equiv c\equiv d \pmod{2} \right\}.$$ If $\bigo^{\times}$ denotes the group of units of $\bigo$, we have that $\bigo^{\times}$ is generated (as a multiplicative group) by $i$, $j$, and $\frac{1+i+j+k}{2}$, and that $|\bigo^{\times}| = 24$. 

For any positive integer $n$, we define the $n$-th Hecke operator $T_n : C(S^2) \rightarrow C(S^2)$ by $$ (T_n f)(z) := \frac{1}{|\bigo^{\times}|}\sum_{\substack{\gamma \in \bigo \\ \mathrm{Nrd}(\gamma) = n}} f(\gamma \cdot z).$$ When $n$ is a power of two, the $n$-th Hecke operator degenerates to either $T_1$ or $T_2$ depending on the $2$-adic valuation of $n$. Hence, the Hecke operators of interest are the ones corresponding to odd values of $n$. We restrict to these for the present discussion.

Since the Laplacian is rotation-invariant, all the Hecke operators commute with $\Delta$. A classic result from number theory is that Hecke operators are self-adjoint and commute with each other. Thus, we can diagonalize further a basis of $H_{\ell}$ with respect to the Hecke operators to obtain an orthonormal basis of spherical harmonics which are eigenfunctions of all the Hecke operators. We call these Hecke eigenfunctions (or arithmetic spherical harmonics).

From the definition of $T_n$, we can see that if $f\in C(S^2)$, then $T_n f$ is invariant under the action of $\bigo^{\times}$, that is, $$(T_n f)(\gamma \cdot z) = (T_n f)(z) \text{ for all } \gamma \in \bigo^{\times}, \: z\in S^2.$$ Hence, Hecke eigenfunctions are elements of the subspace $H_{\ell}^{\bigo^{\times}} \subset H_{\ell}$ consisting of functions invariant under the action of $\bigo^{\times}$ and are naturally defined on the quotient space $M=\bigo^{\times} \backslash S^2$. In fact, the spectral resolution of $\Delta$ in $L^2(\bigo^{\times} \backslash S^2)$ is in terms of these Hecke eigenfunctions. This is the setup where we focus from now on.

Now, let $D<0$ be a fundamental discriminant with $D\not\equiv 1 \pmod{8}$. Consider the set $$ \widehat{\mathcal{E}}(n) := \left\{ \left( \frac{x}{\sqrt{n}}, \frac{y}{\sqrt{n}}, \frac{z}{\sqrt{n}} \right) \in \R^3 : (x,y,z) \in \mathcal{E}(n) \right\} \subset S^2,$$ where 
\begin{equation*}
    n = n_D := \begin{cases}
        -D/4 & \text{ if } 4\mid D, \\
        -D & \text{ otherwise}
    \end{cases}    
\end{equation*}
and $$ \mathcal{E}(n) := \left\{ (x,y,z)\in \Z^3 : x^2 + y^2 + z^2 = n \right\}.$$ The elements of $\widehat{\mathcal{E}}(n)$ are called the \textit{CM points} of $S^2$ (associated to $D$). We require $D\not\equiv 1 \pmod{8}$ so that $\mathcal{E}(n) \neq \varnothing$. Since the actions of $i$, $j$, and $\frac{1+i+j+k}{2}$ on the sphere are given by
\begin{align*}
    i &: (x,y,z) \mapsto (x,-y,-z), \\
    j &: (x,y,z) \mapsto (-x,y,-z), \\
    \frac{1+i+j+k}{2} &: (x,y,z) \mapsto (z,x,y),
\end{align*}
and all of them preserve $\widehat{\mathcal{E}}(n)$, we obtain an action of $\bigo^{\times}$ on $\widehat{\mathcal{E}}(n)$. 

For a point $\xi \in S^2$, we let $\gamma_\xi$ be the geodesic with north pole at $\xi$ and let $\overline{\gamma_{\xi}}$ be its projection to $\bigo^{\times} \backslash S^2$, which is a geodesic of $\bigo^{\times} \backslash S^2$. We define $$ \mathcal{C}_{D} := \bigcup_{\xi \in \bigo^{\times} \backslash \widehat{\mathcal{E}}(n)} \overline{\gamma_{\xi}} \subset \bigo^{\times} \backslash S^2.$$ Our main result is the following.

\begin{theorem}\label{MainTheorem}
Let $\ell\geq 1$ be an integer and let $\psi\in H_{\ell}^{\bigo^{\times}}$ be an $L^2$-normalized Hecke eigenfunction. Let $\mathcal{C}_{D}$ be as defined above. Then, for all $\varepsilon>0$, we have 
\begin{equation}
    ||\psi|_{\mathcal{C}_{D}}||_{L^2} \ll_{D,\varepsilon} \ell^{\varepsilon}. \label{MainTheorem-eq1}
\end{equation}
\end{theorem}

\begin{remarks}$ $
\begin{enumerate}
    \item[(1)] In broad terms, the strategy of the proof is to express the square of this restricted $L^2$-norm in terms of a weighted moment of $L$-functions, which is then bounded using techniques from analytic number theory. This strategy was first pioneered by Ghosh, Reznikov, and Sarnak in \cite{GhoshReznikovSarnak}, who studied the same problem for Hecke--Maass cusp forms on the modular surface $\mathrm{SL}_{2}(\mathbb{Z}) \backslash \mathbb{H}$ restricted to the vertical geodesic joining $0$ and $i\infty$; here the connection to $L$-functions is via the Plancherel theorem for the Mellin transform and the theory of the Hecke integral. For closed geodesics on the modular surface, this problem was studied by Ali in \cite{DanaMScThesis}, where now the connection to $L$-functions is via Waldspurger's formula for $\mathrm{GL}_2$. In the setting studied in this paper, we instead use Waldspurger's formula for inner forms of $\mathrm{GL}_2$.
    \item[(2)] Our result crucially relies on the fact that the Ramanujan conjecture is known for the Hecke eigenvalues of arithmetic spherical harmonics. This is a consequence of the Jacquet--Langlands correspondence (see Section \ref{SpherHarmonics}) and Deligne's proof of the Ramanujan conjecture for holomorphic modular forms (see \cite{Deligne-WeilConj}). Without this result, our bound would instead take the weaker form $$||\psi|_{\mathcal{C}_D}||_{L^2} \ll_{D,\varepsilon} \ell^{\frac{7}{64} + \varepsilon},$$ which matches Ali's bound \eqref{AbouAli-bound}.
    \item[(3)] The method of proof of Theorem \ref{MainTheorem} only applies to geodesics whose north pole is a CM point. One could hope to use Theorem \ref{MainTheorem} and the classic result of Linnik\footnote{Linnik showed in \cite{Linnik} that CM points equidistribute on the sphere if we restrict to the ones that correspond to fundamental discriminants $D$ that satisfy $n_D \equiv \pm 1 \pmod{5}$. The result including all fundamental discriminants was later showed independently by Duke and Schulze-Pillot in \cite{DukeSchulzePillot} and Golubeva and Fomenko in \cite{GolubevaFomenko}.} that CM points equidistribute on the sphere to obtain a similar bound for general geodesics. The issue with this approach is that the implicit constant in \eqref{MainTheorem-eq1} depends polynomially on $|D|$, so that if we fix the Hecke eigenfunction and let $|D|\rightarrow \infty$, the upper bound in \eqref{MainTheorem-eq1} goes to $\infty$ as well.
    \item[(4)] If we assume the generalized Lindel{\"o}f hypothesis, then we would get the same bound as in Theorem \ref{MainTheorem}. Hence, from this perspective, our result is essentially sharp. 
    \item[(5)] In the case of the modular surface $\mathrm{SL}_{2}(\mathbb{Z}) \backslash \mathbb{H}$, a well-known conjecture of Berry \cite{Berry1977} and of Hejhal and Rackner \cite{HejhalRackner} posits that Laplacian eigenfunctions ought to behave like random waves (i.e.~Gaussian random functions) in the large Laplacian eigenvalue limit. This (loosely-defined) conjecture is known as the random wave conjecture; see \cite{SarnakAQC} for an extended discussion on this conjecture and arithmetic quantum chaos. If one assumes the analogue of the random wave conjecture for Laplacian eigenfunctions on $\bigo^{\times} \backslash S^2$, then the correct upper bound in \eqref{MainTheorem-eq1} would be $\ll 1$. However, due to the arithmeticity of Hecke eigenfunctions, it is expected that Hecke eigenfunctions can achieve large values at CM points, yielding to large values of these restricted $L^2$-norms on arithmetic geodesics that are not predicted by the random wave conjecture. In fact, this already happens for the modular surface; see \cite{Michels} and \cite{Milicevic}. This phenomenon for $\bigo^{\times} \backslash S^2$ will be explored further in future work.
\end{enumerate}
\end{remarks}

The structure of this paper is the following. In Section \ref{Background}, we recall some of the important facts that we shall need for the proof of Theorem \ref{MainTheorem}. In Section \ref{Proof}, we prove Theorem \ref{MainTheorem}.

\section{Background}\label{Background}

\subsection{Spherical harmonics and automorphic representations of \texorpdfstring{$B^{\times}(\A)$}{B(A)}}\label{SpherHarmonics}

The space $L^2(S^2)$ decomposes as $$ L^2(S^2) = \bigoplus_{\ell=0}^{\infty} H_\ell,$$ where, as before, $H_\ell$ is the eigenspace of the Laplacian eigenvalue $\lambda_{\ell} = \ell(\ell+1)$, which is explicitly described as the space of harmonic homogeneous polynomials in three variables of degree $\ell$ restricted to $S^2$. Moreover, we know the dimension and a basis of $H_\ell$. For this, recall that for $n\geq 1$ an integer, the $n$-th Legendre polynomial is defined by $$P_n(x) := \frac{1}{2^n n!} \frac{\mathrm{d}^n}{\mathrm{d}x^n}(x^2-1)^n.$$ For each integer $m$ with $1\leq m\leq n$, define $$ P_{n}^{m}(x) := (-1)^m (1-x^2)^{m/2} \frac{\mathrm{d}^m}{\mathrm{d}x^m} P_n(x).$$ The functions $P_{n}^{m}(x)$ are commonly called associated Legendre polynomials. Using these polynomials, we have the following basis of $H_\ell$, which follows from \cite[\S 54--56]{Hobson} and \cite[\S VII.9]{MacRobert}.

\begin{proposition}\label{dimBasisHk}
Let $\ell\geq 0$. Then $\dim_{\C} {H_\ell} = 2\ell+1$ and $$ \{P_\ell(\cos(\phi))\} \cup \{P_{\ell}^{n}(\cos(\phi))\cos(n\theta) \: |\: 1\leq n\leq \ell\} \cup \{P_{\ell}^{n}(\cos(\phi))\sin(n\theta) \: |\: 1\leq n\leq \ell\}$$ is an $L^2(S^2)$-orthogonal basis for $H_\ell$, where $$(x,y,z)=(\sin(\phi) \cos(\theta), \sin(\phi) \sin(\theta), \cos(\phi))$$ for $\theta \in [0,2\pi)$ and $\phi \in [0,\pi]$.
\end{proposition}

As mentioned before, our focus is on the surface $M=\bigo^{\times} \backslash S^2$. For this space, we have the decomposition $$ L^2(\bigo^{\times} \backslash S^2) = \bigoplus_{\ell=0}^{\infty} H_{\ell}^{\bigo^{\times}}$$ where, as defined above, $H_{\ell}^{\bigo^{\times}}$ is the subspace of functions $\psi\in H_{\ell}$ such that $\psi(\gamma \cdot z) = \psi(z)$ for all $\gamma \in \bigo^{\times}$ and $z \in S^2$.

Now, take $\psi \in H_{\ell}^{\bigo^{\times}}$. We can lift $\psi$ to an eigenfunction $\widehat{\psi}$ of the Casimir operator\footnote{It is a distinguished element of the center of the universal enveloping algebra of the (complex) Lie algebra associated to $B^{\times}(\R)$. Intuitively, the Casimir operator is an analogue of the Laplacian operator for Riemannian manifolds in the context of Lie algebras. For more on this operator, see \cite[\S 10.2]{HallLieStuff}.} on $B^{\times}(\R)$ via the formula $$ \widehat{\psi}(x) := \psi(x\cdot k),$$ where $k=(0,0,1)$ is the north pole of $S^2$. In fact, if $K^{\infty} = \mathrm{Stab}(k) \cap \mathrm{Nrd}^{-1}(1)$, the function $\widehat{\psi}$ is defined on $\bigo^{\times} \backslash B^{\times}(\R)/K^{\infty}$ since $\psi$ is taken to be left-invariant under $\bigo^{\times}$.

For any ring $R$, we define $$B(R) := \langle 1,i,j,k \: | \: i^2=j^2=k^2=-1, \: ij=-ji=k\rangle_R,$$ so that $B(R) = B(\Q) \otimes_{\Q} R$ for any $\Q$-algebra $R$. We denote by $B^{\times}(R)$ the group of units of $B(R)$. Then $B^{\times}$ is an affine algebraic group over $\Q$. For $\alpha = a+bi+cj+dk \in B(R)$, let $\nrd{\alpha} := a^2+b^2+c^2+d^2 \in R$ and $\trd{\alpha} := 2a \in R$ be its (reduced) norm and (reduced) trace, respectively.

For a global field $K$, let $\A_{K}$ be the ring of $K$-adèles and $\A_{K,\text{fin}}$ the ring of finite $K$-adèles; we denote by $\A := \A_{\Q}$ and $\A_{\text{fin}} := \A_{\Q,\text{fin}}$. Since $\bigo$ has class number $1$, we have by \cite[eq. (2.20)]{AutRepGetz} the homeomorphism 
\begin{align}
    L : \bigo^{\times} \backslash B^{\times}(\R) &\longrightarrow B^{\times}(\Q) \backslash B^{\times}(\A) / \widehat{\bigo}^{\times} \label{homeo2.20Bcross} \\
    \bigo^{\times} x_{\infty} &\longmapsto B^{\times}(\Q) (x_{\infty},1,1,1,\dots ) \widehat{\bigo}^{\times}, \nonumber
\end{align}
where $\widehat{\bigo}$ is the closure of $\bigo$ in $B(\A_{\text{fin}})$. Thus, we can lift $\widehat{\psi}$ to an automorphic form on $B^{\times}(\A)$ given by $L_{*} \widehat{\psi} = \widehat{\psi} \circ L^{-1}$. If $\psi$ is a Hecke eigenfunction, then $L_{*} \widehat{\psi}$ generates an irreducible cuspidal automorphic representation of $B^{\times}(\A)$ of trivial central character, denoted by $\rho_{\psi}$. Moreover, $L_{*} \widehat{\psi}$ is the unique $K^{\infty} \widehat{\bigo}^{\times}$-invariant vector in $\rho_{\psi}$. Now, by Flath's tensor product theorem (see \cite[Thm. 6.3.4]{AutRepGetz}), we have the tensor product decomposition 
\begin{equation}
    \rho_{\psi} = \bigotimes_{p\leq \infty} {\vphantom{\sum}}' \rho_{p}. \label{Flath}
\end{equation}
The local representations $\rho_p$ are ramified only at $2$ and $\infty$. The representation $\rho_{\infty}$ has as a model the $\ell$-th irreducible representation of $\R^{\times} \backslash B^{\times}(\R) \cong \mathrm{SO}_3(\R)$, which is a representation on the space $H_\ell$ (see \cite[Thm. 17.12]{QuantumTheoryforMath}).

Hecke eigenfunctions have a strong connection with modular forms. Indeed, let $S_{\ell}^{\text{new}}(\Gamma_0(2))$ denote the space of holomorphic newforms of weight $\ell$ and level $2$. By Eichler \cite{EichlerHeckeOp}, we have a map
\begin{align}
    H_{\ell}^{\bigo^{\times}} &\longrightarrow S_{2\ell+2}^{\text{new}}(\Gamma_0(2)) \label{LiftSpherHarmToModForms} \\
    \psi &\longmapsto f_{\psi} \nonumber
\end{align}
given by $$ f_{\psi} (z) := \frac{1}{24}\sum_{\gamma \in \bigo} \left\langle \psi(\xi), \psi\left(\gamma \xi \overline{\gamma} \right) \right\rangle_{L^2(S^2)} e^{2\pi i \nrd{\gamma} z},$$ where $\left\langle \psi(\xi), \psi\left(\gamma \xi \overline{\gamma} \right) \right\rangle$ denotes the inner product in $L^2(S^2)$. Eichler showed that $f_{\psi}$ is an eigenform if $\psi$ is a Hecke eigenfunction and, if this is the case and $T_n \psi = \lambda_{\psi}(n) \psi$, then $T_n f_{\psi} = \lambda_{\psi}(n) n^{\ell} f_{\psi}$ for $n$ odd. In fact, the map \eqref{LiftSpherHarmToModForms} is an explicit realization of the Jacquet--Langlands correspondence between $B^{\times}$ and $\mathrm{GL}_{2}$ (see \cite{GlobalJL}, \cite{JLOriginal}) in the sense that, if $\pi_{\psi}$ is the Jacquet--Langlands transfer of $\rho_{\psi}$ to a cuspidal automorphic representation of $\mathrm{GL}_{2}(\mathbb{A})$, then the adèlic lift of $f_{\psi}$ to an adèlic automorphic form on $\mathrm{GL}_{2}(\mathbb{A})$ generates $\pi_{\psi}$. This correspondence tells us that the Hecke eigenspaces in $H_{\ell}^{\bigo^{\times}}$ are one-dimensional.

\subsection{Closed geodesics of the 2-sphere}

As we mentioned before, we will work with geodesics on $\bigo^{\times} \backslash S^2$. Before this, we first recall what the geodesics of $S^2$ are. These are the equators of $S^2$ when we choose a point of $S^2$ to be the ``north pole''.

\begin{theorem}\label{geodesicsS2}
$\mathrm{(}$\cite[Example 2.80(c)]{RiemannianGeometryGHL}$\mathrm{)}$ The nonconstant geodesics of $S^2$ are the ``great circles'', i.e.~the circles which result of the intersection of $S^2$ with any plane of $\R^3$ containing the origin. In particular, all the geodesics of $S^2$ are closed.
\end{theorem}

\begin{figure}[!htb]
    \centering
    \includegraphics[width=0.6\linewidth]{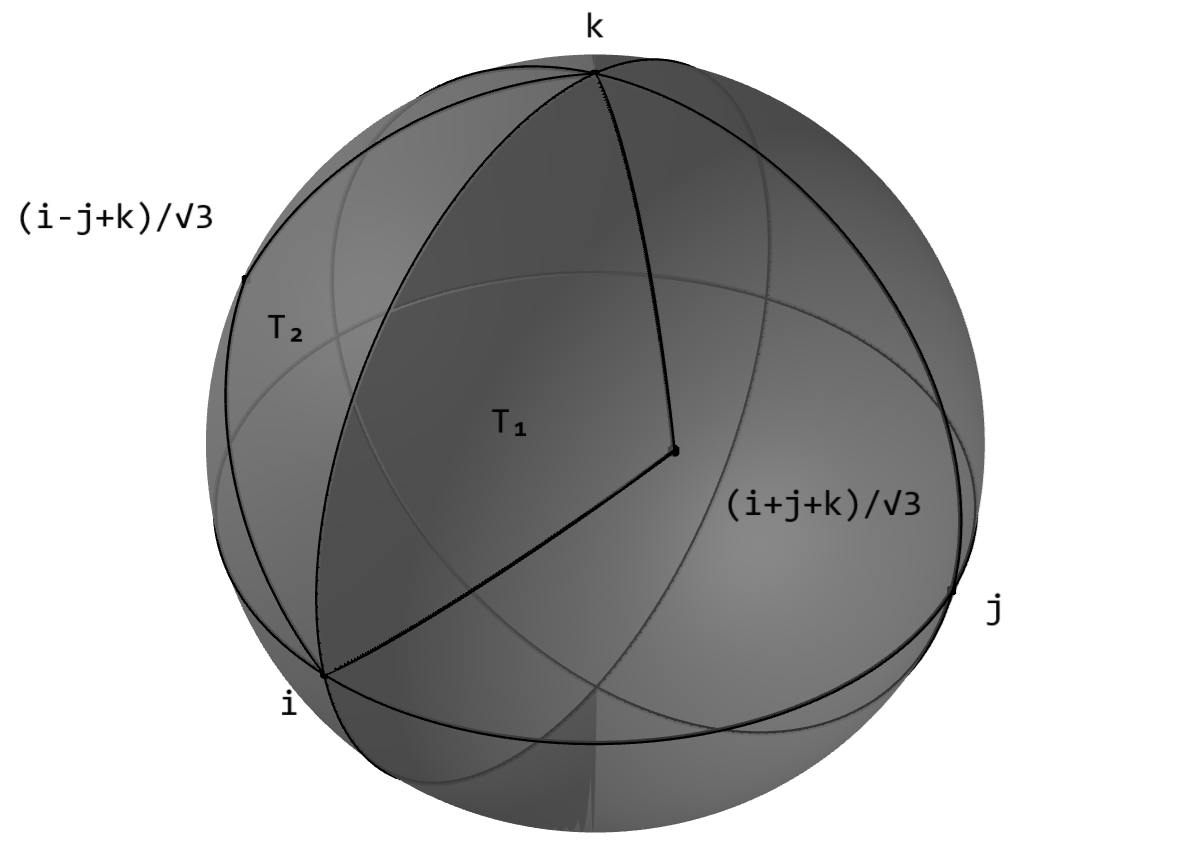}
    \caption{A fundamental domain for the action of $\bigo^{\times}$ on $S^2$ is $T_1 \cup T_2$.}
    \label{MageeSphere}
\end{figure}

The group $\bigo^{\times}$ acts properly on $S^2$ by isometries, but it does not act freely ($-1$ and $1$ act trivially, and any other element of $\bigo^{\times}$ has exactly two fixed points which are antipodal). However, the quotient space $\bigo^{\times} \backslash S^2$ is a smooth manifold and inherits a Riemannian metric from $S^2$. From \cite{Magee}, a fundamental domain for $\bigo^{\times} \backslash S^2$ is given by $T_1 \cup T_2$, where $T_1$ and $T_2$ are the spherical triangles from Figure \ref{MageeSphere}.
In terms of the geodesics of $\bigo^{\times} \backslash S^2$, since the Riemannian metric on $\bigo^{\times} \backslash S^2$ is inherited from the Riemannian metric on $S^2$, using \cite[Thm. 14.8]{DiffGeoTu} we get the following.

\begin{proposition}\label{geodesicsOS2}
The nonconstant geodesics of $\bigo^{\times} \backslash S^2$ are the projections of the geodesics of $S^2$ on $\bigo^{\times} \backslash S^2$.
\end{proposition}

The length of any geodesic of $S^2$ is $2\pi$. The length of its projection to $\bigo^{\times} \backslash S^2$ need not have the same length. In fact, one can show the following result.

\begin{proposition}\label{lengthGeodesics}
Let $\gamma$ be a geodesic of $S^2$ and $\overline{\gamma}$ the corresponding geodesic of $\bigo^{\times} \backslash S^2$. Let $\xi\in S^2$ be a ``north pole'' of $S^2$ with equator $\gamma$, i.e.~a point on $S^2$ such that the geodesic $\gamma$ is the intersection of $S^2$ with the plane $\langle \xi , (x,y,z)\rangle =0$. Denote by $\ell(\gamma)$ and $\ell(\overline{\gamma})$ the lengths of $\gamma$ and $\overline{\gamma}$, respectively, in the corresponding spaces.
\begin{enumerate}[(1)]
    \item If exactly two entries of $\xi$ are zero, then $\ell(\overline{\gamma}) = \frac{1}{4}\ell(\gamma)$.
    \item If exactly one entry of $\xi$ is zero, then $\ell(\overline{\gamma}) = \frac{1}{2}\ell(\gamma)$.
    \item If the absolute values of the entries of $\xi$ are all equal, then $\ell(\overline{\gamma}) = \frac{1}{2} \ell(\gamma)$.
    \item If none of the above occurs, then $\ell(\overline{\gamma}) = \ell(\gamma)$.
\end{enumerate}
\end{proposition}

\subsection{CM modular forms}

There is an important construction of modular forms from a given Hecke character of an imaginary quadratic field $\Q(\sqrt{D})$. It was introduced by Hecke in \cite{HeckeDihedral}. This construction satisfies that, if $\chi$ is the Hecke character and $f_{\chi}$ is the corresponding modular form, then $L(s,f_{\chi}) = L(s,\chi)$, where $L(s,\chi)$ is seen as a degree $1$ $L$-function over $\Q(\sqrt{D})$. In other words, $f_{\chi}$ is the automorphic induction of $\chi$ from $\Q(\sqrt{D})$ to $\mathbb{Q}$.

For our exposition, we follow \cite[\S 12.3]{Iwaniec}. Throughout this section, we let $E=\Q(\sqrt{D})$ be an imaginary quadratic field of discriminant $D<0$ and $\bigo_E$ be its ring of integers.

\begin{definition}\label{Def-dihedralMaass}
Let $\mathfrak{m}$ a modulus of $E$ and $\xi$ a primitive Hecke character modulo $\mathfrak{m}$ for which there exists a non-negative integer $u$ such that $$ \xi((a)) = \left(\frac{a}{|a|}\right)^{u} \text{ if } a\equiv 1 \pmod{\mathfrak{m}}.$$ We define the theta series associated to $\xi$, denoted by $f_{\xi}$, by $$ f_{\xi}(z) := \sum_{\mathfrak{a}} \xi(\mathfrak{a}) N_{E/\Q}(\mathfrak{a})^{u/2} e^{2\pi i z\cdot N_{E/\Q}(\mathfrak{a})},$$ where the sum is taken over all the nonzero integral ideals of $\bigo_E$. The function $f_{\xi}$ is called the \textit{CM modular form} associated to $\xi$.
\end{definition}

\begin{theorem}\label{dihedralMaassCusp}
$\mathrm{(}$\cite[Thm. 12.5]{Iwaniec}$\mathrm{)}$ We have that $f_{\xi}$ is a holomorphic newform of weight $k = u+1$, level $N=|D|\cdot N_{E/\Q}(\mathfrak{m})$, and character $\chi \pmod{N}$, which is given by $$ \chi(n) = \chi_D(n) \xi((n)) \text{ for } n\in \Z,$$ where $\chi_D$ is the primitive quadratic character modulo $-D$ associated to $E$. Moreover, if $u>0$, then $f$ is a cusp form.
\end{theorem}

\subsection{\texorpdfstring{$L$}{L}-functions}

We recall some notation: given $\pi$ an automorphic representation of $\mathrm{GL}_n(\A_F)$, we denote by $L(s,\pi)$ the associated automorphic $L$-function (which is an $L$-function of degree $n$ over $F$), $\gamma(s,\pi)$ its gamma factor, $\epsilon(f\otimes g)$ its epsilon factor, $q(\pi)\in \N$ the conductor of $\pi$, and $$ \mathfrak{q}_{\infty}(s,\pi) := q(\pi) \prod_{i=1}^{d} (3+|s+\kappa_i|)$$ the analytic conductor of $L(s,\pi)$, where the $\kappa_i$ are the local parameters of $\pi$ at infinity.

The first type of $L$-function that we will work with are Rankin--Selberg $L$-functions, possibly defined over a quadratic extension of $\Q$ (see \cite[\S 5.1, 5.11]{IwaniecKowalski}). However, in our specific case, this $L$-function will be equal to a Rankin--Selberg $L$-function defined over $\Q$. To be more precise, using properties of base change, automorphic induction, and Rankin--Selberg $L$-functions, one can show the following.

\begin{proposition}\label{TwistAutoRepEtoF}
Let $E/F$ be a quadratic extension of number fields. Let $\pi$ be a cuspidal automorphic representation of $\mathrm{GL}_{2}(\A_F)$ and $\Omega:E^{\times} \backslash \A_{E}^{\times} \rightarrow \C^{\times}$ a Hecke character. Let $\pi_E = \mathcal{BC}_{E/F}(\pi)$ be the base change of $\pi$ to an automorphic representation of $\mathrm{GL}_2(\A_{E})$ and $\mathcal{AI}_{E/F}(\Omega)$ be the automorphic induction of $\Omega$ to an automorphic representation of $\mathrm{GL}_{2}(\A_{F})$. Then 
\begin{equation}
    L(s,\pi_E \otimes \Omega) = L(s,\pi\otimes \mathcal{AI}_{E/F}(\Omega)), \label{TwistAutoRepEtoF-0}
\end{equation}
where the left hand side of \eqref{TwistAutoRepEtoF-0} is a degree $2$ $L$-function over $E$, and the right hand side of \eqref{TwistAutoRepEtoF-0} is a degree $4$ $L$-function over $F$.
\end{proposition}

The second type of $L$-functions that we will deal with are adjoint $L$-functions associated to newforms. The basic facts of these can be found in \cite[\S 5.12]{IwaniecKowalski}. We will mainly need to understand the behavior of $L(s,\mathrm{Ad}(f))$ at $s=1$. For this purpose, we have the following.

\begin{theorem}\label{AdjointBound}
$\mathrm{(}$\cite[Thm. 0.2]{AdjointBound}$\mathrm{)}$ Let $f$ be a holomorphic newform of weight $k$ and level at most $2$. Then $$ L(1,\mathrm{Ad}(f)) \gg_{\varepsilon} k^{-\varepsilon}$$ for any $\varepsilon>0$.
\end{theorem}

\subsection{Approximate functional equation}

As we will see later, the sums that we will estimate will involve central values of Rankin--Selberg $L$-functions associated to self-dual newforms, which can be well approximated by the approximate functional equation. The general facts can be found in \cite[\S 5.2]{IwaniecKowalski} and its errata. For our specific situation, we have the following results.

\begin{theorem}[Approximate functional equation]\label{AppFuncEq1}
Let $f,g$ be self-dual newforms such that the completed $L$-function $\Lambda(s,f\otimes g)$ is entire. Assume that $q(f)$ and $q(g)$ are square-free and coprime. Let $G(u) = e^{u^2}$. We have $$L(1/2, f\otimes g) = (1+\epsilon(f\times g)) \sum_{m\geq 1} \frac{\chi_f(m) \chi_g(m)}{m} \sum_{n\geq 1} \frac{\lambda_f(n) \lambda_g(n)}{\sqrt{n}} V_{1/2} \left( \frac{m^2 n}{\sqrt{q(f\times g)}} \right),$$ where $$V_s(y) := \frac{1}{2\pi i} \int_{\mathrm{Re}(u)=3} \frac{\gamma(s+u,f\otimes g)}{\gamma(s,f\otimes g)} y^{-u} G(u) \: \frac{\mathrm{d}u}{u}$$ and $\epsilon(f\otimes g)$ is the root number of $L(s,f\otimes g)$.
\end{theorem}

\begin{remark}\label{AppFuncEqHoldsForDeven}
The conditions on $q(f)$ and $q(g)$ help us know what the Dirichlet series of $L(s,f\otimes g)$ is. For our purposes, we will use Theorem \ref{AppFuncEq1} when $f$ has level $2$ with trivial nebentypus and $g$ has level $|D|$ with primitive nebentypus, where $D<0$ is a fundamental discriminant. When $D$ is odd, the conditions on $q(f)$ and $q(g)$ are satisfied; when $D$ is even, they are not. However, the theorem still holds since the Dirichlet series of $L(s,f\otimes g)$ is still as illustrated in Theorem \ref{AppFuncEq1} by \cite[Thm. 2.2]{WinnieLiRankinSelberg}.
\end{remark}

In fact, one can show that, for our choice of $G(u)$, the function $V_s(y)$ decays rapidly.

\begin{lemma}\label{DecayVsy}
Fix $A>0$ and $\alpha$. Suppose that $\mathrm{Re}(s+\kappa_j) \geq \alpha >0$ for $j=1,2,3,4$, where $\kappa_j$ are the local parameters of $L(s,f\otimes g)$ at infinity. Then, the function $V_s(y)$ defined in Theorem \ref{AppFuncEq1} satisfies $$ V_s(y) \ll_{\alpha,A} \left(1+\frac{y}{\sqrt{\q_{\infty}(s,f\otimes g)}} \right)^{-A}.$$
\end{lemma}

Thus, in the approximate functional equation of Theorem \ref{AppFuncEq1}, we can ignore all the terms that are large when compared to the analytic conductor $\q_{\infty}(1/2,f\otimes g)$.

\begin{corollary}\label{AppFuncEq2}
Consider the same notation and assumptions from Theorem \ref{AppFuncEq1}. Furthermore, assume that $\epsilon(f\otimes g) = 1$. Then, for any $B,c,\delta,M>0$ such that $M \geq c\cdot \q_{\infty}(1/2,f\otimes g)^{\frac{1}{2}+\delta}$ we have 
\begin{align*}
    &L(1/2, f\otimes g) \\
    &= 2 \sum_{m=1}^{M} \frac{\chi_f(m) \chi_g(m)}{m} \sum_{n=1}^{M} \frac{\lambda_f(n) \lambda_g(n)}{\sqrt{n}} V_{1/2} \left( \frac{m^2 n}{\sqrt{q(f\times g)}} \right) + O_{c,B,\delta,q}\left(\q_{\infty}(1/2,f\otimes g)^{-B}\right),
\end{align*} 
where $V_s(y)$ is the function defined in Theorem \ref{AppFuncEq1} and the implied constant depends only on $c,B,\delta$, and $q:=q(f\otimes g)$.
\end{corollary}

\subsection{Dyadic partition of unity}

It is common in analytic number theory that a sum over $\Z$ needs to be analyzed over smaller intervals, since it becomes easier to estimate. For this, we use partitions of unity.

\begin{proposition}\label{PartitionUnity}
There exists a smooth function $W(x)$ supported on $[1/2,2]$ such that $$ \sum_{k\geq 0} W\left( \frac{x}{2^k} \right) = 1 $$ for any $x\geq 1$.
\end{proposition}

\subsection{Lipschitz principle}

As described by Davenport in \cite{Davenport-LipschitzPrin}, the idea of the Lipschitz principle is to count lattice points inside a sufficiently nice region. For our purposes, we will need the region to be a parallelogram. For this, we have the following result.

\begin{proposition}\label{Lipschitz-Dana}
$\mathrm{(}$\cite[Cor. 3.5]{DanaMScThesis}$\mathrm{)}$ Let $\Lambda$ be a lattice in $\R^2$ with (ordered) basis $B = \{\alpha_1, \alpha_2\}$. Let $v_1, v_2\in \R^2$ be two linearly independent vectors. Let $\mathcal{P}:= \{av_1 + bv_2 : a,b\in [0,1]\}$. Then $$ |\mathcal{P}\cap \Lambda| \ll_{B} (1+||v_1||)(1+||v_2||),$$ where the implied constant depends only on the basis $B$ of $\Lambda$.
\end{proposition}

\section{Proof of Theorem \ref{MainTheorem}}\label{Proof}

\subsection{Setup}

Let $\psi$ be an $L^2$-normalized Hecke eigenfunction. Let $D<0$ be a negative fundamental discriminant and assume that $D\not\equiv 1 \pmod{8}$. As mentioned in Section \ref{Introduction}, the elements of $\mathcal{C}_D$ are geodesics on $\bigo^{\times} \backslash S^2$ corresponding to the CM points of $S^2$ associated to $D$. One can count the number of geodesics in $\mathcal{C}_D$. For this, we start with a classical result of Gauss.

\begin{theorem}[Gauss]\label{GaussThm}
Assume $n>4$, square-free, and $n\not\equiv 7\pmod{8}$. Let $r_3(n) := \# \mathcal{E}(n) = \# \widehat{\mathcal{E}}(n)$. Then
\begin{equation*}
    r_3(n) = \begin{cases}
        12h_D & \text{ if } n \equiv 1,2 \pmod{4}, \\
        24h_D & \text{ if } n \equiv 3 \pmod{8},
    \end{cases}
\end{equation*}
where $h_D$ is the ideal class number of the imaginary quadratic field $E:=\Q(\sqrt{D})$.
\end{theorem}

In \cite{VenkovGaussThm}, Venkov gave a different proof of Theorem \ref{GaussThm} using the arithmetic of the Hurwitz integers $\bigo$.\footnote{A similar result can be obtained for the number of representations of positive integers by certain ternary quadratic forms using other quaternion algebras with a maximal order of class number $1$; see \cite{ShemanskeActionCLEonEn}.} Later, in \cite{RehmActionCLEonEn}, Rehm reproved Gauss' Theorem using Venkov's method, but written in a more modern framework. They showed that there is a free group action $$\Delta_{D} : \mathrm{Cl}(E) \times \bigo^{\times} \backslash \widehat{\mathcal{E}}(n) \longrightarrow \bigo^{\times} \backslash \widehat{\mathcal{E}}(n),$$ from where they obtained that
\begin{equation*}
    |\bigo^{\times} \backslash \widehat{\mathcal{E}}(n)| = \left\{ \begin{array}{ll}
        h_D & \text{ if } n\equiv 1,2 \pmod{4} \\
        2h_D & \text{ if } n\equiv 3 \pmod{8}.
    \end{array} \right.
\end{equation*}
Let $$c_D := \frac{|\bigo^{\times} \backslash \widehat{\mathcal{E}}(n)|}{h_D} = \left\{ \begin{array}{ll}
    1 & \text{ if } D\equiv 0 \pmod{4}, \\
    2 & \text{ if } D\equiv 5 \pmod{8}.
\end{array} \right..$$
Thus, by fixing an element in each orbit of the action of $\mathrm{Cl}(E)$ on $\bigo^{\times} \backslash \widehat{\mathcal{E}}(n)$ and letting $\mathrm{Cl}(E)$ act on them, we obtain $c_D$ maps $\{\tau_1, \dots, \tau_{c_D}\}$ from $\mathrm{Cl}(E)$ to $\bigo^{\times}\backslash \widehat{\mathcal{E}}(n)$, where any two maps have disjoint image and the union of the images of all the maps is $\bigo^{\times} \backslash \widehat{\mathcal{E}}(n)$. Hence, we obtain a bijective map
\begin{align*}
    \tau : \Z/c_D \Z \times \mathrm{Cl}(E) &\longrightarrow \bigo^{\times} \backslash \widehat{\mathcal{E}}(n) \\
    ([b],[\a]) &\longmapsto \tau_{b'}([\a]),
\end{align*}
where $b' \equiv b \pmod{c_D}$ and $1\leq b' \leq c_D$. For $[b] \in \Z/c_D \Z$ and $[\a] \in \mathrm{Cl}(E)$, denote by $\mu_{[b],[\a]}$ a lift of $\tau([b],[\a])$ to a point in $\widehat{\mathcal{E}}(n)$ and let $\gamma_{[b],[\a]}$ be the corresponding geodesic of $S^2$.

\subsection{\texorpdfstring{$L^2$}{L2}-restriction norm formula}

Let $\gamma:=\gamma_{\mu}$ be a geodesic of $S^2$ associated to $\mu \in \widehat{\mathcal{E}}(n)$ and let $\overline{\gamma}$ be its projection to $\bigo^{\times} \backslash S^2$. Since $\psi$ is an $L^2$-normalized Hecke eigenfunction defined on $\bigo^{\times} \backslash S^2$, by the proof of Proposition \ref{lengthGeodesics}, we can see that
\begin{equation}
    ||\psi|_{\overline{\gamma}}||_{L^2} = c_{\gamma} ||\psi|_{\gamma}||_{L^2}, \label{eq-1}
\end{equation}
where $c_{\gamma}\in \{1,1/2,1/4\}$. 
Therefore, we want to bound 
\begin{equation*}
    ||\psi|_{\overline{\gamma}} ||_{L^2}^2 = c_{\gamma} \int_{\gamma} |\psi(x)|^2 \: \mathrm{d}s.
\end{equation*} 

Let $g:\overline{\gamma} \rightarrow \C$ be any continuous function and consider its lift to $\gamma$, which we also denote by $g$ by abuse of notation. Let $\kappa_{\gamma} \in B^{\times}(\R)$ be such that $\kappa_{\gamma} \cdot \mu = i$. A parametrization of $\gamma$ is given by $\theta \mapsto \kappa_{\gamma}^{-1} \cdot e^{i \theta} k$ with $\theta \in [0,2\pi)$, where $k$ is the north pole of $S^2$. Then, we have
\begin{equation}
    \int_{\gamma} g(x) \: \mathrm{d}s = \int_{0}^{2\pi} g\left(\kappa_{\gamma}^{-1} \cdot e^{i \theta} k \right) \: \mathrm{d}\theta. \label{eq1}
\end{equation}
Let 
\begin{align*}
    \Phi : [0,2\pi) & \rightarrow \gamma_i \\
    \theta & \mapsto e^{i \theta}k.
\end{align*}
Then, the right hand side of \eqref{eq1} is the integral of the pullback $\Phi^{*} g( \kappa_{\gamma}^{-1} \cdot )$ along $[0,2\pi)$. Canonically identifying $[0,2\pi)$ with $\R/2\pi\Z$ with the Haar probability, $\Phi$ is then a homeomorphism between $\gamma_i$ and a compact abelian group. Moreover, by \eqref{eq1}, this homeomorphism satisfies 
\begin{equation}
    \langle g_1, g_2 \rangle = \langle \Phi^{*} g_1, \Phi^{*} g_2 \rangle, \label{eq2}
\end{equation}
where $g_1,g_2 : \gamma \rightarrow \C$ are any continuous functions and the inner products are in the corresponding $L^2$-spaces. 

Let $$\widetilde{\mathcal{C}_D} := \bigcup_{[b]\in \Z/c_D \Z} \bigcup_{[\a] \in \mathrm{Cl}(E)} \gamma_{[b],[\a]},$$ which is a lift of $\mathcal{C}_D$ to $S^2$. We can parametrize $\widetilde{\mathcal{C}_{D}}$ as
\begin{align*}
    \Xi: \R/(2\pi \Z) \times \Z/c_D \Z \times \mathrm{Cl}(E) &\longrightarrow \widetilde{\mathcal{C}_D} \\
    ([\theta], [b], [\a]) & \longmapsto \kappa_{[b],[\a]}^{-1} \cdot e^{i\theta}k,
\end{align*}
where $\kappa_{[b],[\a]} = \kappa_{\gamma_{[b],[\a]}}$. Let $G_D := \R/(2\pi \Z) \times \Z/c_D \Z \times \mathrm{Cl}(E)$. We endow $\Z/c_D \Z$ and $\mathrm{Cl}(E)$ with the discrete topologies, so that the counting measures become the Haar measures. Consequently, by \eqref{eq2}, the map $\Xi$ is a homeomorphism between $G_D$ and $\widetilde{\mathcal{C}_D}$ compatible with integration; that is, for any integrable function $g:\widetilde{\mathcal{C}_D} \rightarrow \C$, we have $$ \int_{\widetilde{\mathcal{C}_D}} g(x) \: \mathrm{d}s = \sum_{b=1}^{c_D} \sum_{[\a] \in \mathrm{Cl}(E)} \int_{\gamma_{[b],[\a]}} g(x) \: \mathrm{d}s = \int_{G_D} \Xi^{*}g([\theta], [b], [\a]) \: \mathrm{d}([\theta], [b], [\a]).$$ For the inner products of the $L^2$-spaces, we have
\begin{equation}
    \langle g_1, g_2\rangle  = \langle \Xi^{*} g_1, \Xi^{*} g_2\rangle. \label{InnerProdEquiv}
\end{equation}

We know that $G_D$ is a compact abelian group, so its set of characters, $\widehat{G_D}$, is an orthogonal basis of $L^2(G_D)$ (see \cite[Prop. 5.23]{FollandHarmonic}). Thus, we can express the $L^2$-norm of any $g\in L^2(\mathcal{C}_{D})$ in terms of the characters of $G_D$. After identifying $\R/2\pi\Z$ with $S^1$ via $\Phi$, the group of characters of $\R/2\pi\Z$ is $$ \left\{ \chi_n : x \mapsto e^{i nx} : n\in \Z \right\}.$$ By \cite[Cor. 4.8]{FollandHarmonic}, we have that $$ \widehat{G_D} = \left\{ \chi_{\lambda,\rho,n} : \lambda \in \widehat{\Z/c_D \Z}, \rho \in \widehat{\mathrm{Cl}(E)}, n\in \Z \right\},$$ where
\begin{align*}
    \chi_{\lambda,\rho,n} : G_D &\longrightarrow S^1\\
    ([\theta], [b], [\a]) &\longmapsto e^{in\theta} \lambda([b]) \rho([\a]).
\end{align*}

\begin{proposition}
For any $g\in L^2(\mathcal{C}_{D})$, we have
\begin{equation}
    ||g||_{L^2}^{2} \asymp \frac{1}{2\pi c_D h_D} \sum_{\lambda \in \widehat{\Z/c_D \Z}} \sum_{\rho \in \widehat{\mathrm{Cl}(E)}} \sum_{n\in \Z} I(g,\lambda, \rho, n), \label{Planch-eq1}
\end{equation}
where
\begin{equation}
    I(g,\lambda, \rho, n) = \left| \sum_{[b] \in \Z/c_D \Z} \lambda([b])^{-1} \sum_{[\a]\in \mathrm{Cl}(E)} \rho([\a])^{-1} \int_{0}^{2\pi} g\left( \kappa_{[b], [\a]}^{-1} \cdot e^{i \theta} k \right) e^{-i n\theta} \: \mathrm{d}\theta \right|^2. \label{Planch-eq2}
\end{equation}
\end{proposition}

\begin{proof}
From \eqref{eq-1}, \eqref{InnerProdEquiv}, and Parseval's identity for $L^2(\widetilde{\mathcal{C}_{D}})$, since $c_{\gamma_{[b], [\a]}} \asymp 1$, we have that $$ ||g||_{L^2}^{2} \asymp \sum_{\lambda \in \widehat{\Z/c_D \Z}} \sum_{\rho \in \widehat{\mathrm{Cl}(E)}} \sum_{n\in \Z} \frac{|\langle \Xi^{*}g, \chi_{\lambda, \rho,n}\rangle|^2}{||\chi_{\lambda,\rho,n}||_{L^2}^{2}}.$$ The result follows from the fact that the $L^2$-norm of each character of $G_D$ is $\sqrt{2\pi c_D h_D}$.
\end{proof}

In general, the expression for $||g||_{L^2}^{2}$ in \eqref{Planch-eq1} is an infinite sum. However, in our particular case, this infinite sum reduces to a finite sum. To show this claim, we do the following reductions.

\begin{remark}\label{PsiRealValued}
Notice that, by properties of rotations and spherical harmonics, $\overline{\psi}$ is also a Hecke eigenfunction with the same Laplace eigenvalue and same Hecke eigenvalues as $\psi$. Since we know the joint eigenspaces are one-dimensional, it must happen that $\overline{\psi} = c \psi$ for some constant $c\in \C$ with $|c|=1$. This is only possible if $\mathrm{Im}(\psi)$ is contained in a line passing through the origin, since the equation $\overline{z}=cz$ (for $z\in \C$) represents a line passing through the origin. Thus, there exists a constant $c'\in \C^{\times}$ with $|c'|=1$ such that $c' \psi$ takes only real values.
\end{remark}

\begin{lemma}
Let $\psi$ be a spherical harmonic of degree $\ell$. Then, for all $\lambda\in \widehat{\Z/c_D \Z}$, $\rho\in \widehat{\mathrm{Cl}(E)}$, and $n\in \Z$, we have $$I(\psi,\overline{\lambda}, \overline{\rho},-n) = I(\psi,\lambda, \rho, n).$$
\end{lemma}

\begin{proof}
The result follows from taking the complex conjugate of the expression inside the absolute value in the definition of $I(g,\overline{\lambda}, \overline{\rho}, -n)$ and using Remark \ref{PsiRealValued}.
\end{proof}

\begin{lemma}\label{IntEq0forLargen}
Let $\psi$ be a spherical harmonic of degree $\ell$. Then, for all $\lambda\in \widehat{\Z/c_D \Z}$, $\rho\in \widehat{\mathrm{Cl}(E)}$, and $n\in \Z$ such that $n>\ell$, we have $I(\psi,\lambda, \rho, n)=0$.
\end{lemma}

\begin{proof}
For each integer $n\geq 1$, let $$P_n(x) := \frac{1}{2^{n} n!} \frac{\mathrm{d}^n}{\mathrm{d}x^n} (x^2-1)^n$$ be the $n$-th Legendre polynomial and, for each $0\leq m\leq n$, let $$ P_{n}^{m}(x) := (-1)^m (1-x^2)^{m/2} \frac{\mathrm{d}^{m}}{\mathrm{d}x^{m}} P_n(x)$$ be the associated Legendre polynomials. Let $h(\xi) := \psi(\kappa_{[b],[\a]}^{-1} \cdot \xi)$. Since the Laplacian operator of $S^2$ is rotation-invariant, we know that $h \in H_\ell$. By Proposition \ref{dimBasisHk}, $h$ is a linear combination of $P_{\ell}(\cos(\phi))$, $P_{\ell}^{m}(\cos(\phi)) \cos(m\theta')$, and $P_{\ell}^{m}(\cos(\phi)) \sin(m\theta')$, where $1\leq m\leq \ell$ and $\xi = (x,y,z) = (\cos(\theta') \sin(\phi), \sin(\theta') \sin(\phi), \cos(\phi))$ is in spherical coordinates. Hence, to prove that the desired integral vanishes, it suffices to show it vanishes for $h(\xi)$ equal to some element of this basis.

Notice that the integral in the definition of $I(\psi,\lambda, \rho, n)$ is the integral of $h(\xi)$ along the equator of $S^2$ with north pole $(1,0,0)$. Then, in spherical coordinates, the coordinate $\theta'$ remains constant equal to $-\pi/2$ or to $\pi/2$. Thus, to prove the result, we only need to show it for $P_\ell(\cos(\phi))$ and for $P_{\ell}^{m}(\cos(\phi))$, where $1\leq m\leq \ell$. Since $e^{i\theta}k = -j\sin(\theta) + k\cos(\theta)$, then $\phi=\theta$ for the point $e^{i\theta}k$ on $S^2$, so that we want to show, for $n>\ell$ and for $1\leq m\leq \ell$, that $$ \int_{0}^{2\pi} P_{\ell}(\cos(\theta)) e^{-in\theta} \: \mathrm{d}\theta = \int_{0}^{2\pi} P_{\ell}^{m}(\cos(\theta)) e^{-in\theta} \: \mathrm{d}\theta=0.$$ However, $P_\ell(\cos(\theta))$ and $P_{\ell}^{m}(\cos(\theta))$ are polynomials in $\sin(\theta)$ and $\cos(\theta)$ of degree $\ell$, so that they can be expressed as linear combinations of the functions $$e^{-i\ell \theta}, e^{-i(\ell-1)\theta}, \dots, e^{-i\theta}, 1, e^{i\theta}, \dots, e^{i(\ell-1)\theta}, e^{i\ell \theta}.$$ Since each one of these functions is orthogonal to $e^{in\theta}$ in $L^2([0,2\pi])$ for $n>\ell$, the result follows.
\end{proof}

Putting together the previous two results, we obtain the following.

\begin{corollary}\label{L2FormulaFinite}
Let $\psi$ be an $L^2$-normalized Hecke eigenfunction of degree $\ell$. Then $$ ||\psi|_{\mathcal{C}_{D}}||_{L^2}^{2} \asymp \frac{1}{c_D h_D} \sum_{\lambda \in \widehat{\Z/c_D \Z}} \sum_{\rho \in \widehat{\mathrm{Cl}(E)}} \sum_{n=0}^{\ell} I(g,\lambda, \rho, n).$$
\end{corollary}

\begin{remark}\label{ParityObs}
In fact, using symmetry, we can see that $I(\psi,\lambda, \rho, n)=0$ if $n$ and $\ell$ have different parity. This fact is used to simplify computations in the case $D=-4$, but for the general case it is not needed.
\end{remark}

\subsection{Adèlization}\label{Adelization}

Next, we proceed to interpret the expression on the right hand side of \eqref{Planch-eq2} as a period integral. First, we need certain constructions. As before, let $E=\Q(\sqrt{D}) = \Q(\sqrt{-n_D})$. Let $\bigo_E$ be the ring of integers of $E$ and $\widehat{\bigo_E}$ be the closure of $\bigo_E$ in $\A_{E,\text{fin}}$. Note that $$\bigo_E = \left\{ \begin{array}{ll}
    \Z[\sqrt{-n_D}] & \text{ if } D\equiv 0 \pmod{4}, \\
    \Z\left[ \frac{1+\sqrt{-n_D}}{2} \right] & \text{ if } D \equiv 1 \pmod{4}
\end{array} \right.$$ and $$\bigo_{E}^{\times} = \left\{ \begin{array}{ll}
    \{\pm 1\} & \text{ if } D<-4, \\
    \{\pm 1, \pm i\} & \text{ if } D=-4, \\
    \{\pm 1, \pm \omega, \pm \omega^2\} & \text{ if } D=-3,
\end{array} \right.$$ where $\omega = (-1+\sqrt{-3})/2 = e^{2\pi i/3}$ is a primitive cube root of unity. For a prime ideal $\p$ of $\bigo_E$, let $E_{\p}$ be the completion of $E$ at $\p$ and $\bigo_{\p}$ be the valuation ring of $E_{\p}$. We have $\widehat{\bigo_{E}} = \prod_{\p} \bigo_{\p}$, so that $\widehat{\bigo_{E}}^{\times} = \prod_{\p} \bigo_{\p}^{\times}$.

Let $Q(n_D) := \{\mu \in \bigo : \mu^2 = -n_D\}$ and let $\mu \in Q(n_D)$. Then $\nrd{\mu}^2 = \nrd{\mu^2} = n_D^2$, so that $\mu \overline{\mu} = \nrd{\mu} = n_D$. Since $B(\R)$ is a division ring and $\mu^2=-n_D$, it follows that $\overline{\mu} = -\mu$, and so $\trd{\mu}=0$. From here, it follows that there is a bijection
\begin{align*}
    \mathcal{E}(n_D) &\longrightarrow Q(n_D) \\ 
    (x,y,z) &\longmapsto xi + yj + zk.
\end{align*}
The group $\bigo^{\times}$ acts on $Q(n_D)$ by conjugation and the actions of $\bigo^{\times}$ on $\mathcal{E}(n_D)$ and $Q(n_D)$ commute with the bijection above. Hence, we get a correspondence between orbits in $\bigo^{\times} \backslash \mathcal{E}(n_D)$ and orbits in $\bigo^{\times} \backslash Q(n_D)$.

Moreover, assume we have a $\Q$-algebra homomorphism $\iota : E \rightarrow B(\Q)$. It follows that $$\iota(\sqrt{-n_D})^2 = \iota(-n_D)=-n_D.$$ Of course, the image of $\sqrt{-n_D}$ determines $\iota$. Moreover, if $\iota, \iota' : E\rightarrow B(\Q)$ are $\Q$-algebra homomorphisms that satisfy that the images of $\sqrt{-n_D}$ under $\iota$ and $\iota'$ are in the same orbit of $\bigo^{\times} \backslash Q(n_D)$, then $\iota'$ and $\iota$ are related by an inner automorphism of $B(\Q)$ induced by an element of $\bigo^{\times}$. Thus, we obtain a correspondence between $\bigo^{\times} \backslash \mathcal{E}(n_D)$ and $\Q$-algebra homomorphisms $E \rightarrow B(\Q)$ with $\iota(\sqrt{-n_D})\in Q(n_D)$ related under an inner automorphism by an element of $\bigo^{\times}$.

\begin{lemma}\label{iotaOptEmb}
For any $\Q$-algebra homomorphism $\iota : E \rightarrow B(\Q)$ with $\iota(\sqrt{-n_D})\in Q(n_D)$, we have that $\iota(E) \cap \bigo = \iota(\bigo_E)$; i.e.~$\iota$ is an optimal embedding.
\end{lemma}

\begin{proof}
Let $\mu = \iota(\sqrt{-n_D})$ and suppose $\alpha = a+b\mu \in \bigo$ for some $a,b\in \Q$. Since $\trd{\mu}=0$, it follows that $\trd{\alpha}=a$, and so $a \in \frac{1}{2} \Z$. Then $2a\in \Z$, so that $2b\mu$ has integral coordinates. Since $n_D$ is square-free and the coordinates of $\mu$ are integral, it follows that the coordinates of $\mu$ are pairwise coprime, from where it follows that $2b\in \Z$. Moreover, the same argument shows that if $a\in \Z$, then $b\in \Z$. For $D\equiv 1\pmod{4}$, the result follows from here. So, assume $D\equiv 0\pmod{4}$, so that $-n_D\equiv 2,3\pmod{4}$, i.e.~$n_D\equiv 1,2\pmod{4}$. Now, assume that $a\notin \Z$, so that $a=a'/2$ with $a'$ an odd integer. It follows that $b=b'/2$ with $b'$ an odd integer. In this case, we must have that $\alpha \in \bigo - B(\Z)$, so that all the entries of $\alpha$ must be of the form $m/2$ with $m$ an odd integer. Let $\mu = xi+yj+zk$. Since $n_D\equiv 1,2\pmod{4}$ and $\mu^2=-n$, then $x^2+y^2+z^2=n_D$, so that at least one of $x,y,z$ must be even. This implies that $2\alpha = a' + b'xi + b'yj + b'zk$ where at least one coordinate is even and at least one coordinate is odd, which contradicts our assumptions. It follows that necessarily $a,b\in \Z$, so that $\alpha \in \iota(\bigo_E)$.
\end{proof}

Let $\mu \in Q(n_D)$ and consider the embedding $\iota_{\mu}: E \hookrightarrow B(\Q)$ associated to $\mu$, i.e.~the embedding $a+b\sqrt{-n_D} \mapsto a+b\mu$. Note that, by tensoring with $\R$, $\iota_{\mu}$ extends to an embedding $\widetilde{\iota}_{\mu} : \C \hookrightarrow B(\R)$. Since for any field extension $L/K$, the $K$-algebra $L$ is flat ($L$ is a free $K$-module, and so the property follows by \cite[\S 10.5, Cor. 42]{DummitFoote}), we know that the resulting tensor product of maps $\iota_{\mu} \otimes \mathrm{id}_{\Q_p} : E\otimes_{\Q} \Q_p \rightarrow B(\Q) \otimes_{\Q} \Q_p = B(\Q_p)$ is injective. However, by \cite[\S I.10, pg. 57]{CasselsFrohlich}, we know that 
\begin{equation}
    E\otimes_{\Q} \Q_p \cong \bigoplus_{\ell=1}^{g} E_{p,\ell}, \label{NumberFieldTensorpAdic}
\end{equation}
algebraically and topologically, where $E_{p,1}, \dots, E_{p,g}$ are all the possible completions of $E$ with respect to a prime ideal of $\bigo_E$ lying above the prime $p$. Thus, for each prime $p$, we get an embedding $$ \iota_{\mu,p} : \bigoplus_{\ell=1}^{g} E_{p,\ell} \hookrightarrow B(\Q_p).$$ From the embedding $\widetilde{\iota}_{\mu}$ and the embeddings $\iota_{\mu,p}$, we obtain an embedding $$\widehat{\iota}_{\mu} : \C \times \prod_{\mathfrak{p}} E_{\mathfrak{p}} \hookrightarrow B(\R) \times \prod_{p} B(\Q_p),$$ which restricts to an embedding of (topological) rings 
\begin{equation*}
    \widehat{\iota}_{\mu} : \A_{E} \hookrightarrow B(\A).
\end{equation*} 

Indeed, to show this claim, we just need to show that if $\alpha \in \A_{E,\text{fin}}$, then $\widehat{\iota}_{\mu}(\alpha) \in B(\A_{\text{fin}})$. By the lemma in \cite[\S I.12, pg. 61]{CasselsFrohlich}, the claim reduces to verifying that $\iota_{\mu}(\bigo_E) \subset \bigo$, which is true by Lemma \ref{iotaOptEmb}. 
Moreover, considering $E$ as a (discrete) subring of $\A_{E}$ and $B(\Q)$ as a (discrete) subring of $B(\A)$, we can see that $\widehat{\iota}_{\mu}(E) \subset B(\Q)$ since $\widehat{\iota}_{\mu}$ restricts to $\iota_{\mu}$ on $E$. Consequently, we obtain an embedding
\begin{equation*}
    E^{\times} \backslash \A_{E}^{\times} \hookrightarrow B^{\times}(\Q) \backslash B^{\times}(\A).
\end{equation*}

\subsection{Hecke characters}

We need to construct unitary characters on $\A_{E}^{\times}$ that will allow us to interpret the expressions in \eqref{Planch-eq2} as period integrals. From now on, for the ease of exposition, we assume $D<-4$; the cases $D=-3,-4$ can be handled in a similar way with some modifications. First, by \cite[eq. (2.20)]{AutRepGetz}, we have that 
\begin{equation}
    \bigsqcup_{[\a] \in \mathrm{Cl}(E)} \Gamma_{[\a]} \backslash \C^{\times} \cong E^{\times} \backslash \A_{E}^{\times} / \widehat{\bigo_E}^{\times}, \label{AdelicIso1}
\end{equation}
where $\Gamma_{[\a]} = \bigo_E^{\times}$ for all $[\a] \in \mathrm{Cl}(E)$ since $\A_{E}^{\times}$ is abelian. Using the explicit descriptions of the isomorphism in \eqref{AdelicIso1}, we can see that the diagram
\begin{equation}
\begin{tikzcd}
	{\displaystyle \bigsqcup_{[\a] \in \mathrm{Cl}(E)} \Gamma_{[\a]} \backslash \mathbb{C}^{\times}} & {E^{\times} \backslash \mathbb{A}_{E}^{\times} / \widehat{\mathcal{O}_E}^{\times}} \\
	{\Z^{\times} \backslash \R^{\times}} & {\mathbb{Q}^{\times} \backslash \mathbb{A}_{\mathbb{Q}}^{\times} / \widehat{\Z}^{\times}},
	\arrow["\cong", from=1-1, to=1-2]
	\arrow[hook, from=2-1, to=1-1]
	\arrow["\cong"', from=2-1, to=2-2]
	\arrow[hook, from=2-2, to=1-2]
\end{tikzcd} \label{DiagEmbeddings}
\end{equation}
is commutative, where the left vertical map is inclusion into the component corresponding to the trivial class in $\mathrm{Cl}(E)$.

For each $n\in \Z$, consider the homomorphism $\omega_n : \C^{\times} \rightarrow S^1$ defined by $$ \omega_n(z) := \left( z/\overline{z} \right)^n.$$ Since $D<-4$, we have that $\bigo_{E}^{\times} = \{ \pm 1\}$. We know that $\omega_n|_{\bigo_{E}^{\times}}$ is trivial, so that $\omega_n$ can be seen as a character of $\bigo_E^{\times} \backslash \C^{\times}$. Using the isomorphism \eqref{AdelicIso1}, where $\omega_n$ acts on each component of the left side of \eqref{AdelicIso1}, we can lift $\omega_n$ to a unitary Hecke character $\Omega_n : \A_{E}^{\times} \rightarrow \C^{\times}$ that is trivial on $\widehat{\bigo_E}^{\times}$ and has infinity type $\omega_n$. Moreover, the diagram \eqref{DiagEmbeddings} implies that $\Omega_n$ is trivial on $\A^{\times}$ since $\omega_n$ is trivial on $\R^{\times}$.

\begin{remark}\label{D34HeckeCharacters}
For the cases $D=-3$ and $D=-4$, the homomorphism $\omega_n$ is trivial on $\bigo_{E}^{\times}$ only for certain values of $n$ (for $n$ even if $D=-4$ and for $n$ a multiple of $3$ if $D=-3$), so that the same construction does not work. To fix this issue, we consider a modulus $\mathfrak{f}$ of $E$ such that the units in $\bigo_{E}^{\times}$ are congruent to $1 \pmod{\mathfrak{f}}$, and then proceed to construct a Hecke character of conductor dividing $\mathfrak{f}$ following \cite[\S 6, pp. 9--10]{HeckeCharJerryShurman}. In both cases, we can take $\mathfrak{f} = (3)$. We shall see how this change affects the proof later.
\end{remark}

\subsection{Period integrals}

Let $\pi^B:=\rho_{\psi}$ be the (cuspidal) automorphic representation of $B^{\times}(\A)$ generated by $L_{*} \widehat{\psi}$ as described in Section \ref{SpherHarmonics}. Fix $n\in \Z$, $[b] \in \Z/c_D \Z$, and $[\a]\in \mathrm{Cl}(E)$. The adèlic object that we are interested in is the normalized period integral $|P^{B}(\phi)|^2 / \langle \phi, \phi \rangle$. Here, for $\phi \in \pi^{B}$, we define
\begin{equation}
    P^{B}(\phi) := \int_{E^{\times} \A^{\times} \backslash \A_{E}^{\times}} \phi(\widehat{\iota}_{\mu_{[b],[\a]}}(t)) \overline{\Omega_{n}}(t) \: \mathrm{d}t, \label{PBPeriodInt}
\end{equation}
with the measure $\mathrm{d}t$ being the Tamagawa measure on $\A_{E}^{\times}$, which is normalized so that $E^{\times} \A^{\times} \backslash \A_{E}^{\times}$ has volume $2 \Lambda(1, \chi_{D})$, where $\chi_{D}$ is the primitive quadratic character modulo $-D$ associated to $E$, and 
\begin{equation}
    \langle \phi, \phi\rangle := \int_{Z(\A) B^{\times}(\Q) \backslash B^{\times}(\A)} |\phi(g)|^2 \: \mathrm{d}^{\times}g, \label{InnProd}
\end{equation}
where $\mathrm{d}^{\times}g$ is the Tamagawa measure, which is normalized so that $Z(\A) B^{\times}(\Q) \backslash B^{\times}(\A)$ has volume $2$. Note that the period integral \eqref{PBPeriodInt} is well-defined by the properties of $\phi$ and the construction of $\Omega_n$.

There are two reasons to be interested in the period integral \eqref{PBPeriodInt}. One of them is that we can express the square of its absolute value in terms of $L$-functions via Waldspurger's formula. The second is that it is closely related to the integral on the right hand side of \eqref{Planch-eq2}. 

Indeed, let $\phi_0$ be the unique spherical vector of $\pi^B = \rho_{\psi}$ (i.e.~the adèlic lift of $\psi$ to $B^{\times}(\A)$). Let $\kappa_{[b],[\a]}'\in B^{\times}(\R)$ be such that $$\nrd{\kappa_{[b],[\a]}'} = 1, \quad \kappa_{[b],[\a]}'\cdot k = \kappa_{[b],[\a]}^{-1}\cdot k, \quad \text{ and } \quad (\kappa_{[b],[\a]}')^{-1} \cdot k \in \gamma_i,$$ where $\gamma_i$ is the geodesic of $S^2$ with north pole $i$. Such $\kappa_{[b],[\a]}' \in B^{\times}(\R)$ exists since there exists $\mu' \in \gamma_i$ such that $d(k,\mu') = d(\kappa_{[b],[\a]}^{-1} \cdot k,k)$, where $d$ denotes the distance function in $S^2$, so that there exists $\kappa \in B^{\times}(\R)$ such that\footnote{This observation follows from the fact that $S^2$ is a \textit{two-point homogeneous space}; see \cite[\S 5.2, Thm. II]{TwoPointHomSpaces}.} $\kappa \cdot k = \kappa_{[b],[\a]}^{-1}\cdot k$ and $\kappa \cdot \mu' = k$; then $\kappa_{[b],[\a]}' = \kappa / \sqrt{\nrd{\kappa}}$ satisfies the desired properties. Define $\phi_{[b], [\a]}(g) := \phi_0\left( g\kappa_{[b], [\a]}' \right)$, which is in $\pi^{B}$ since $\kappa_{[b],[\a]}'$ is in the maximal compact subgroup of $B^{\times}(\R)$.

We know that both $\phi_0$ and $\Omega_n$ are invariant under $\widehat{\bigo_E}^{\times}$, hence $\phi_{[b], [\a]}$ as well. Since the normalization of the Tamagawa measure $\mathrm{d}t$ is such that the measure of $\widehat{\bigo_E}^{\times}$ is $1$, we have $$ P^B(\phi_{[b], [\a]}) = \int_{E^{\times} \A^{\times} \backslash \A_{E}^{\times} / \widehat{\bigo_E}^{\times}} \phi_{[b], [\a]}(\widehat{\iota}_{\mu_{[b],[\a]}}(t)) \overline{\Omega_{n}}(t) \: \mathrm{d}t.$$ However, by \eqref{AdelicIso1}, we then obtain that $$ E^{\times} \A^{\times} \backslash \A_{E}^{\times} / \widehat{\bigo_E}^{\times} \cong \bigsqcup_{[\mathfrak{c}] \in \mathrm{Cl}(E)} \Gamma_{[\mathfrak{c}]} \backslash S^1,$$ algebraically and topologically\footnote{The left hand side has the quotient topology induced from the idèlic topology, whereas the right hand side has the disjoint union topology and each component has the quotient topology induced from the usual topology on $S^1$.}. 
In fact, we can make the homeomorphism explicit. Let
\begin{align*}
    \mathrm{Cl}(E) &\longrightarrow E^{\times} \backslash \A_{E,\text{fin}}^{\times} / \widehat{\bigo_E}^{\times} \\
    [\mathfrak{c}] &\longmapsto E^{\times} g_{[\mathfrak{c}]} \widehat{\bigo_{E}}^{\times}
\end{align*}
be a bijection. Then
\begin{align*}
    \bigsqcup_{[\mathfrak{c}] \in \mathrm{Cl}(E)} \Gamma_{[\mathfrak{c}]} \backslash S^1 &\longrightarrow E^{\times} \A^{\times} \backslash \A_{E}^{\times} / \widehat{\bigo_E}^{\times} \\
    \Gamma_{[\mathfrak{c}]} z &\longmapsto E^{\times} \A^{\times} g_{[\mathfrak{c}]} z \widehat{\bigo_E}^{\times}
\end{align*}
is a homeomorphism. Since the Tamagawa measure is a Haar measure, we can supress the $g_{[\mathfrak{c}]}$ in the integral. Moreover, by the definition of the embedding $\widehat{\iota}_{[b],[\a]}$ and the fact that $\Gamma_{[\mathfrak{c}]}=\bigo_{E}^{\times} = \{\pm 1\}$, we obtain that $$ \left\{ \kappa_{[b],[\a]}^{-1} e^{i\theta} \kappa_{[b],[\a]} : \theta \in [0,\pi] \right\}$$ is a fundamental domain for the component corresponding to $g_{[\mathfrak{c}]}=1$ of $\widehat{\iota}_{[b],[\a]}\left(E^{\times} \A^{\times} \backslash \A_{E}^{\times} / \widehat{\bigo_E}^{\times}\right)$ in $B^{\times}(\A)$. Consequently, by definition of $\kappa_{[b],[\a]}'$,
\begin{align*}
    P^{B}(\phi_{[b], [\a]}) &= \frac{2\Lambda(1,\chi_{D})}{\pi} \sum_{[\mathfrak{c}] \in \mathrm{Cl}(E)} \int_{0}^{\pi} \phi_{[b], [\a]}\left( \kappa_{[b],[\a]}^{-1} e^{i\theta} \kappa_{[b],[\a]} \right) \overline{\Omega_{n}}\left( e^{i\theta} \right) \: \mathrm{d}\theta \\
    &= \frac{2\Lambda(1,\chi_{D}) h_D}{\pi} \int_{0}^{\pi} \psi\left( \left(\kappa_{[b],[\a]}^{-1} e^{i\theta}\right) \cdot (\kappa_{[b],[\a]} \kappa_{[b],[\a]}' \cdot k) \right) e^{-2in\theta} \: \mathrm{d}\theta \\
    &= \frac{2\Lambda(1,\chi_{D}) h_D}{\pi} \int_{0}^{\pi} \psi\left( \left(\kappa_{[b],[\a]}^{-1} e^{i\theta}\right) \cdot k \right) e^{-2in\theta} \: \mathrm{d}\theta \\
    &= \frac{2\Lambda(1,\chi_{D}) h_D}{\pi} \int_{0}^{\pi} \psi\left( \kappa_{[b],[\a]}^{-1} \cdot e^{2i\theta} k \right) e^{-2in\theta} \: \mathrm{d}\theta \\
    &= \frac{\Lambda(1,\chi_{D}) h_D}{\pi} \int_{0}^{2\pi} \psi\left( \kappa_{[b],[\a]}^{-1} \cdot e^{i\theta} k \right) e^{-in\theta} \: \mathrm{d}\theta
\end{align*}
which is precisely the integral in \eqref{Planch-eq2} with $g=\psi$.

On the other hand, for the same choice of $\phi$ (i.e.~$\phi = \phi_{[b], [\a]}$), we can explicitly compute \eqref{InnProd}. Indeed, similarly as in the previous computation, we have that $\phi$ is invariant under $\widehat{\bigo}^{\times}$, where $\bigo$ is the ring of Hurwitz integers in $B(\Q)$ and $\widehat{\bigo}$ is its closure in $B(\A_{\text{fin}})$. The Tamagawa measure is constructed in such a way that the volume of $\widehat{\bigo}^{\times}$ in $B^{\times}(\A_{\text{fin}})$ is $1$. It follows that $$ \langle \phi , \phi \rangle = \int_{Z(\A) B^{\times}(\Q) \backslash B^{\times}(\A) / \widehat{\bigo}^{\times}} |\phi(g)|^{2} \: \mathrm{d}^{\times}g.$$ However, by the homeomorphism \eqref{homeo2.20Bcross}, the construction of $\phi$ from $\psi$, the fact that $\kappa_{[b], [\a]}'$ is a rotation of $S^2$, and considering the normalization of volumes, we obtain that $$ \langle \phi, \phi\rangle = \frac{1}{2\pi} \int_{S^2} |\psi(z)|^{2} \: \mathrm{d}\sigma(z) = \frac{1}{2\pi}$$ since $\psi$ is $L^2(S^2)$-normalized.

\subsection{Formulae relating period integrals and special values of \texorpdfstring{$L$}{L}-functions}

Let $\pi$ be the cuspidal automorphic representation of $\GL_2(\A)$ corresponding to $f_{\psi}$, so that $\pi$ has trivial central character. Let $\pi^{B}$ be the automorphic representation of $B^{\times}(\A)$ corresponding to $\pi$ under the Jacquet--Langlands correspondence, so that $\pi^B = \rho_{\psi}$. Note that, by our constructions, $\pi$ and $\Omega = \Omega_n$ have disjoint ramification. Let $\pi_E = \mathcal{BC}_{E/\Q}(\pi)$ be the base change of $\pi$ to an automorphic representation of $\GL_2(\A_{E})$, which is cuspidal unless $\pi$ is dihedral (see \cite[Lemma 11.3, Prop. 11.4]{LanglandsBaseChangeGL2}). As we mentioned in the previous section, in \cite{WaldspurgerFormula} (specifically Proposition 7), Waldspurger proved a formula relating $|P^{B}(\phi)|^2$ and $L(1/2, \pi_E\otimes \Omega)$. However, the factors in this formula are not as explicit as one would want them. Jacquet and Chen, in \cite[\S 6, Thm. 2]{JacChen}, gave a similar formula with more explicit factors that works in the case of $\pi$ not being dihedral. Martin and Whitehouse, in the appendix to \cite{MartinWhitehouse}, then showed an analogous formula to Jacquet and Chen's for the dihedral case. To use these results, we need to specify a certain distribution which will help us compute the local factors in the formula. We follow \cite{Magee} and \cite{JacChen}.

Let $C_{c}^{\infty}(B^{\times}(\A))$ be the space of smooth, compactly supported, complex-valued functions on $B^{\times}(\A)$. The distribution is given by 
\begin{equation}
    J_{\pi^B} : C_{c}^{\infty}(B^{\times}(\A)) \rightarrow \C, \quad f\mapsto \sum_{\varphi} \int [\pi^{B}(f)\varphi](t) \overline{\Omega}(t) \: \mathrm{d}t \overline{\int \varphi(t) \overline{\Omega}(t) \: \mathrm{d}t}, \label{DefDistribution}
\end{equation}
where the sum is over an orthonormal basis of $\pi^{B}$ and $$ \pi^{B}(f)\varphi = \int_{B^{\times}(\A)} f(g) [\pi^B(g)\varphi]\: \mathrm{d}g.$$

Let $S_0$ be the finite set of places of $\Q$ containing the infinite place and all places where $\pi^B$ or $\Omega$ are ramified. Let $\A^{S_0}$ denote the ring of adèles over $\Q$ with the entries at the places in $S_0$ missing. In our case, since we are assuming that $D<-4$, we have $S_0=\{2,\infty\}$. 
Denote by $S$ the set of places of $E$ lying over $S_0$. For any $f=\left(\prod_{v\in S_0} f_v\right)f^{S_0}$, where $f^{S_0}$ is the characteristic function of a maximal compact subgroup of $B^{\times}(\A^{S_0})$, the distribution $J_{\pi^{B}}$ factorizes at $f$. We take $f$ in such a way that $\pi^{B}(f)$ is the orthogonal projection onto the span of the $L^2(B^{\times}(\A))$-normalized spherical vector $\phi^{*} = \phi_{[b], [\a]}^{*}$ in $\pi^{B}$, where $\phi_{[b], [\a]}^{*} = \phi_{[b], [\a]} / \sqrt{\langle \phi_{[b], [\a]}, \phi_{[b], [\a]} \rangle}$.
Then all the terms in the formula \eqref{DefDistribution} vanish except the one corresponding to $\varphi = \phi^{*}$, in which case we have $\pi^{B}(f) \phi^{*} = \phi^{*}$. This analysis and the analysis from the previous section imply that
\begin{equation}
    J_{\pi^{B}}(f) = |P^{B}(\phi^{*})|^2 = \frac{|P^{B}(\phi_{[b], [\a]})|^2}{\langle \phi_{[b], [\a]}, \phi_{[b], [\a]} \rangle}. \label{Dist+PeriodInt}
\end{equation}

On the other hand, by \cite[\S 6, Thm. 2]{JacChen}, we have that
\begin{equation}
    J_{\pi^{B}}(f) = \prod_{v\in S_0} \tilde{J}_{\pi^B}(f_v) \times \frac{1}{2}\prod_{\substack{v\in S_0 \\ v \text{ inert}}} \varepsilon(1,\eta_{v},\psi_{v}) 2L(0,\eta_{v}) \times \frac{L_{S_0}(1,\eta) L^{S}(1/2, \pi_E \otimes \Omega_n)}{L^{S_0}(1,\mathrm{Ad}(\pi))}, \label{FactJacChen}
\end{equation}
where $\eta$ is the quadratic Hecke character of $\A^{\times}$ associated to $E$, $L^{S}(s,\pi)$ denotes the $L$-function of $\pi$ with the local factors at places in $S$ missing, and $L_{S}(s,\pi)$ denotes the product of the local factors of the $L$-function of $\pi$ at the places in $S$. For the definitions of the undefined terms in the previous formula, see \cite{JacChen}. The important properties of them is that they are positive and depend only on the ramification of $\pi^B$ and of $E$. The same is true for $L_{S_0}(1,\eta)$. The terms $\tilde{J}_{\pi^{B}}(f_{v})$ are local factors of the distribution that need to be calculated case by case, as described in \cite{JacChen}. 

We use formula \eqref{FactJacChen} to estimate the $L^2$-restricted norm of $\psi$. First, we want to do two things: compute the local factors $\tilde{J}_{\pi^B}(f_v)$ and include the missing local factors of the $L$-function corresponding to finite places. We do these in the following two sections.

\begin{remark}
For $D=-3$ and $D=-4$, we have $S_0 = \{2,3,\infty\}$. In both cases, besides the analysis shown in the following two sections, we also need to compute the corresponding local factors at the prime $3$.
\end{remark}

\subsection{Local factors of the distribution}

We now compute the local factors of the distribution on the right hand side of \eqref{FactJacChen}. We first concentrate on $\tilde{J}_{\pi^B}(f_2)$. Under the isomorphism \eqref{Flath}, the vector $\phi^{*}$ corresponds to a pure tensor. Since $\phi^{*}$ is $\widehat{\bigo}^{\times}$-invariant, then the local component of $\phi^{*}$ at $2$ must be invariant under the compact open subgroup $\bigo_{2}^{\times}$ of $B^{\times}(\Q_2)$, where $\bigo_2$ is the closure of $\bigo$ in $B(\Q_2)$ and is a maximal order in $B(\Q_2)$. By \cite[Prop. 2.6]{GrossPrasad1991}, there is a unique vector (up to scalar multiplication) in $\rho_2$ which satisfies this property, and is precisely the Gross--Prasad test vector for linear forms. By the definition of $f$, we have that $\rho_2(f_2)$ is the orthogonal projection onto the span of the local component of $\phi^{*}$ at $2$. With this setup, the factor $\tilde{J}_{\pi^B}(f_2)$ is computed in \cite[\S 2.2, 2.2.1]{MartinWhitehouse} for our case ($B$ non-split at $2$, $\pi$ ramified at $2$), where they obtain that $\tilde{J}_{\pi^B}(f_2)=1$.

We now concentrate on computing $\tilde{J}_{\pi^B}(f_{\infty})$. We follow the proof of \cite[Lemma 8]{Magee}. Now, we assume $|n|\leq \ell$, since the integral $I(\psi,\lambda, \rho, n)$ vanishes if $|n|\geq \ell$. By the discussion in Section \ref{SpherHarmonics}, we can pick as a model for $\rho_{\infty}$ the space $H_{\ell}$ of homogeneous harmonic polynomials of degree $\ell$ on $\R^3$ restricted to the $2$-sphere, which is precisely the $\ell$-th irreducible representation of $\R^{\times} \backslash B^{\times}(\R) \cong \mathrm{SO}_3(\R)$. 

Let $Y_{\ell}^{n}$ be the associated Legendre polynomial $P_{\ell}^{n}$ normalized so that $$||Y_{\ell}^{n}(\langle \cdot, k\rangle )||_{L^2(S^2)}=1,$$ where we use the convention $P_{\ell}^{0} = P_{\ell}$. Let $$ e_1(v) := Y_{\ell}^{n}(\langle v,i\rangle ) e^{in\theta(v)} \: \text{ and } \: e_2(v) := Y_{\ell}^{n}(\langle v,i\rangle ) e^{-in\theta(v)},$$ where $\theta(v)$ is an angle around the axis of $i$ such that $\theta(k)=0$ and $\theta(j)=\pi/2$. These unit vectors are in $H_\ell$ by Proposition \ref{dimBasisHk} and satisfy that $$ \rho_{\infty}(t)e_1 = \overline{\Omega_{n}}(t) e_1 \: \text{ and } \: \rho_{\infty}(t) e_2 = \Omega_n(t) e_2$$ for $t \in \R^{\times} \backslash \C^{\times}$. Indeed, recall that $[\rho_{\infty}(t)e_1](v) = e_1(tvt^{-1})$ for $t\in B^{\times}(\R)$. By a direct computation, we can verify that if $t\in \C^{\times}$, then $\langle tvt^{-1},i\rangle  = \langle v,i\rangle$. Moreover, we can see that if $t=re^{i\theta_t}$ and $\theta(v)=\theta$, then $\frac{v-\langle v,i\rangle i}{||v-\langle v,i\rangle i||} =  ke^{i\theta}$, so that $$ t\left(ke^{i\theta}\right)t^{-1} = re^{i\theta_t} \cdot ke^{i\theta} \cdot r^{-1} e^{-i\theta_t} = ke^{i(\theta-2\theta_t)},$$ from where it follows that $\theta(tvt^{-1}) = \theta - 2\theta_t$. Consequently, $$ [\rho_{\infty}(t)e_1](v) = e_1(tvt^{-1}) = Y_{\ell}^{n}(\langle v,i\rangle) e^{in(\theta-2\theta_t)} = e^{-2in\theta_t} e_1(v) = \overline{\Omega_{n}}(t) e_1(v),$$ and similarly with $e_2$.

From \cite[\S 5.1, \S 6]{JacChen}, since $e_1$ and $e_2$ are elements of an orthonormal basis of $H_{\ell}$, it follows that $$ \tilde{J}_{\pi^B}(f_{\infty}) = \langle \rho_{\infty}(f_{\infty})e_1, e_2\rangle.$$ Now, we chose $f_{\infty}$ in such a way that $\rho_{\psi}(f)$ is the projection onto $\phi^{*}$, which is the translation by $\kappa_{[b],[\a]}'$ of the unique (up to scalar multiplication) $K^{\infty} \widehat{\bigo}^{\times}$-fixed vector of $\rho_{\psi}$. Hence, $\rho_{\infty}(f_{\infty})$ is precisely the orthogonal projection onto the rotation by $\kappa_{[b],[\a]}'$ of the (unique) $K^{\infty}$-fixed element of $H_\ell$, which is $Y_{\ell}^{0}(\langle k,\cdot \rangle)$. Thus, 
\begin{align*}
    [\rho_{\infty}(f_{\infty})e_1](y) &= \left( \int_{S^2} Y_{\ell}^{0}(\langle k,\kappa_{[b],[\a]}' \cdot x\rangle ) Y_{\ell}^{n}(\langle x,i\rangle ) e^{-in\theta} \: \mathrm{d}x \right) Y_{\ell}^{0}(\langle k,\kappa_{[b],[\a]}' \cdot y\rangle) \\
    &= \left( \int_{S^2} Y_{\ell}^{0}(\langle (\kappa_{[b],[\a]}')^{-1} \cdot k, x\rangle ) Y_{\ell}^{n}(\langle x,i\rangle ) e^{-in\theta} \: \mathrm{d}x \right) Y_{\ell}^{0}(\langle (\kappa_{[b],[\a]}')^{-1} \cdot k,y\rangle)
\end{align*}
Since $\sqrt{2\ell+1}Y_{\ell}^{0}(\langle \cdot, \cdot \rangle)$ is the reproducing kernel of $H_{\ell}$ in $L^2(S^2)$ (see \cite[\S 2.3, Lemma 2.22]{MorimotoSphere}), by the choice of $\kappa_{[b],[\a]}'$ we obtain that 
\begin{align*}
    [\rho_{\infty}(f_{\infty})e_1](y) &= \frac{1}{\sqrt{2\ell+1}} \left(Y_{\ell}^{n}(\langle (\kappa_{[b],[\a]}')^{-1} \cdot k,i \rangle) e^{-in\theta((\kappa_{[b],[\a]}')^{-1}\cdot k)}\right) Y_{\ell}^{0}(\langle (\kappa_{[b],[\a]}')^{-1} \cdot k,y\rangle) \\
    &= \frac{1}{\sqrt{2\ell+1}} Y_{\ell}^{n}(0) Y_{\ell}^{0}(\langle (\kappa_{[b],[\a]}')^{-1} \cdot k,y\rangle) e^{-in\theta((\kappa_{[b],[\a]}')^{-1}\cdot k)}.
\end{align*}
Using the reproducing property again, we finally obtain
\begin{align}
    \tilde{J}_{\pi^B}(f_{\infty}) &= \langle \rho_{\infty}(f_{\infty})e_1, e_{2}\rangle \label{ArchPartNonDepab} \\
    &= \frac{Y_{\ell}^{n}(0)}{\sqrt{2\ell+1}} e^{-in\theta((\kappa_{[b],[\a]}')^{-1}\cdot k)} \int_{S^2} Y_{\ell}^{0}(\langle (\kappa_{[b],[\a]}')^{-1} \cdot k,y\rangle ) Y_{\ell}^{n} (\langle y,i\rangle ) e^{in\theta} \: \mathrm{d}y = \frac{Y_{\ell}^{n}(0)^2}{2\ell+1}. \nonumber
\end{align}
In fact, using \cite[Ch. VII, \S 4, eq. (17)]{MacRobert}, we can determine the value of $Y_{\ell}^{n}(0)$.

\begin{lemma}\label{ValueYkn0}
For $0\leq n\leq \ell$ of the same parity, we have $$ Y_{\ell}^{n}(0) = (-1)^{(\ell+n)/2} \frac{\sqrt{(2\ell+1)(\ell+n+2)(\ell-n+2) C_{(\ell+n)/2} C_{(\ell-n)/2}}}{2^{\ell+2}\sqrt{\pi}},$$ where $C_m := \frac{1}{m+1} \binom{2m}{m}$ is the $m$-th Catalan number. If $0\leq n\leq \ell$ have different parity, then $Y_{\ell}^{n}(0)=0$.
\end{lemma}

As a consequence of this result and Stirling's estimates of the factorial, we obtain the following.

\begin{corollary}\label{EstimateArchimFactor}
For $0\leq n\leq \ell$ of the same parity, we have $$ |Y_{\ell}^{n}(0)| \asymp \frac{\sqrt{\ell}}{(1+\ell+n)^{1/4} (1+\ell-n)^{1/4}}.$$
\end{corollary}

From the observations in the previous section and equations \eqref{Dist+PeriodInt}, \eqref{FactJacChen}, and \eqref{ArchPartNonDepab}, we can see that $|P^{B}(\phi_{[b],[\a]})|^2$ does not depend on $[\a] \in \mathrm{Cl}(E)$ and $[b]\in \Z/c_D \Z$, so that for $\lambda \in \widehat{\mathrm{Cl}(E)}$ and $\rho\in \widehat{\mathrm{Cl}(E)}$ we have $$ I(\psi, \lambda, \rho, n) \leq c_{D}^2 h_{D}^2 |P^{B}(\phi_{[b],[\a]})|^2.$$ Consequently,
\begin{equation}
    I(\psi,\rho,\lambda,n) \ll_{D,E,\mathrm{ram}(\Omega_n),\mathrm{ram}(\pi)} \prod_{v\in S_0} \tilde{J}_{\pi^B}(f_v) \times \frac{L^{S}(1/2, \pi_E \otimes \Omega_n)}{L^{S_0}(1,\mathrm{Ad}(\pi))} = \frac{Y_{\ell}^{n}(0)^2}{2\ell+1} \frac{L^{S}(1/2, \pi_E \otimes \Omega_n)}{L^{S_0}(1,\mathrm{Ad}(\pi))}. \label{PeriodIntToLvalues}
\end{equation}

\subsection{Local nonarchimedean factors of the \texorpdfstring{$L$}{L}-functions}

Next, we would like to add the local factors at $2$ to the $L$-functions in our formulas to be able to use the approximate functional equation later. The background for the following discussion can be found in \cite{NTBTate}, \cite[\S 1]{IntroLanglands}, and \cite[\S 7]{localLanglandsGL2}.

First of all, recall that $\pi$ has conductor $2$ and trivial central character since it corresponds to $f_{\psi}$. Thus, by \cite[Prop. 2.8]{LocalComponents}, we know that $\pi_2$ is the twist of the Steinberg representation by a (unitary) unramified character $\chi$ of $\Q_{p}^{\times}$, i.e.~$\pi_2 \simeq \mathrm{St}\otimes \chi$. Consequently, $\pi_2$ corresponds under the local Langlands correspondence to the Weil--Deligne parameters $(\tau, N)$, where $\tau:W_{\Q_2} \rightarrow \GL_2(\C)$ is a representation of the Weil group such that 
\begin{equation*}
    \tau(\Phi) = \begin{pmatrix}
        \chi(\Phi)2^{-1/2} & 0 \\ 0 & \chi(\Phi) 2^{1/2}
    \end{pmatrix},
\end{equation*}
where $\Phi$ is the geometric Frobenius, and 
\begin{equation}
    N = \begin{pmatrix}
        0 & 1 \\ 0 & 0
    \end{pmatrix} \label{NilpotentEndo}
\end{equation}
is a nilpotent endomorphism of $\C^2$. Now, we have the following.

\begin{lemma}\label{2nonsplit}
The prime $2$ of $\Q$ is non-split in $E$.
\end{lemma}

\begin{proof}
If $D$ is even, then $2$ ramifies in $E$. Now, assume $D$ is odd, so $D\equiv 1 \pmod{4}$. We know that $\bigo_{E} = \Z[(1+\sqrt{D})/2]$, where the minimal polynomial of $(1+\sqrt{D})/2$ is $f(x) = x^2 - x + \frac{1-D}{4}$. Recall that we assumed $n\not\equiv 7 \pmod{8}$, so that $D\not\equiv 1 \pmod{8}$, hence $\frac{1-D}{4}$ is odd. Consequently, $f(x)$ is irreducible over $\Z/2\Z$. From here, it follows by Dedekind's Theorem (see \cite[Thm. 3.41]{milneANT}) that there is only one prime ideal lying above $2$ in $E$ and, since $2$ is unramified, it must be inert.
\end{proof}

\begin{remark}
The previous lemma can also be concluded from the fact that $E$ can be embedded in $B(\Q)$. Indeed, consider an embedding $\iota:E \xhookrightarrow{} B(\Q)$ as described in Section \ref{Adelization}. Tensoring with $\Q_2$, we obtain an embedding $\iota' : E\otimes_{\Q} \Q_2 \xhookrightarrow{} B(\Q) \otimes_{\Q} \Q_2 = B(\Q_2)$. If $2$ is split in $E$, then $E\otimes \Q_2 \cong \Q_2 \oplus \Q_2$ by \eqref{NumberFieldTensorpAdic}. However, $\Q_2 \oplus \Q_2$ contains nontrivial idempotents; thus it cannot embed in $B(\Q_2)$ since the latter is a division algebra. Therefore, $2$ must be non-split in $E$. 
\end{remark}

By Lemma \ref{2nonsplit}, there is only one prime ideal $\p$ lying above $2$ in $E$. We analyze the behavior of the local $L$-function $L_{\p}(1/2, \pi_{E}\otimes \Omega_n)$ depending on the factorization of $2$ in $E$.

First, assume $2$ ramifies in $E$, so that $2\bigo_E = \p^2$. Then $E_{\p} \cong \Q_2$, which implies $W_{E_{\p}} \cong W_{\Q_2}$. It follows that $(\pi_{E}\otimes \Omega_n)_{\p}$ is also a twist of the Steinberg representation, now twisted by $\chi$ and by $\Omega_n$. Thus, the Weil--Deligne parameters of $\pi_{E,\p}\otimes (\Omega_n)_{\p}$ are $(\tau',N)$, where $\tau':W_{E_{\p}} \rightarrow \GL_2(\C)$ is a representation of the Weil group such that $$ \tau'(\Psi) = \begin{pmatrix}
    \chi(\Psi) \alpha_n(\Psi) 2^{-1/2} & 0 \\ 0 & \chi(\Psi) \alpha_n(\Psi) 2^{1/2}
\end{pmatrix},$$ where $\alpha_n:\mathrm{Gal}(E_{\p}^{\mathrm{ab}}/E_{\p}) \rightarrow \C^{\times}$ corresponds to $(\Omega_n)_{\p}:E_{\p}^{\times} \rightarrow \C^{\times}$ via local class field theory (see \cite[Ch. III, Thm. 3.1]{milneCFT}) and $\Psi$ is the geometric Frobenius. From here, we have that $\mathrm{Ker}(N) = \{(z,0) \in \C^2\}$, where $\tau'(\Psi)$ acts by multiplication by $\chi(\Psi) \alpha_n(\Psi) 2^{-1/2}$. It follows that $$ L_{\p}(1/2, \pi_{E}\otimes \Omega_n) = \mathrm{det}\left( 1-\tau'(\Psi)|_{\mathrm{Ker}(N)} 2^{-1/2} \right)^{-1} = \left( 1-\frac{\chi(\Psi) \alpha_n(\Psi)}{2} \right)^{-1}.$$ Since the characters $\chi$ and $\alpha_n$ are unitary, we obtain that $L_{\p}(1/2,\pi_E\otimes \Omega_n) \asymp 1$ uniformly on $n$.

Secondly, assume that $2$ is inert in $E$, so that $2\bigo_E = \p$. In this case, $E_{\p}$ is the unique quadratic unramified extension of $\Q_p$ (by the discussion after Prop. 8.10 in \cite{milneANT}) and we have $W_{E_{\p}}$ is an open subgroup of index $2$ in $W_{\Q_p}$. If $\Phi'$ is the geometric Frobenius of $W_{E_{\p}}$, by the definitions of the Weil group $W_{E_{\p}}$ and the geometric Frobenius, we have $\Phi' = \Phi^2$. Additionally, by the definition of base change, $\pi_{E,\p}$ corresponds under the local Langlands correspondence to the Weil--Deligne parameters $(\mathrm{Res}_{W_{E_{\p}}}(\tau),N)$, where $N$ is as defined in \eqref{NilpotentEndo} and $$ \mathrm{Res}_{W_{E_{\p}}}(\tau)(\Phi') = \tau(\Phi^2) = \begin{pmatrix}
    \chi(\Phi)^2 \cdot 2^{-1} & 0 \\ 0 & \chi(\Phi)^2 \cdot 2
\end{pmatrix}.$$ Hence, $$ L_{\p}(1/2, \pi_E\otimes \Omega_n) = \left( 1-\frac{\chi(\Phi)^2 \alpha_n(\Phi)}{2^{3/2}} \right)^{-1},$$ where again $\alpha_n$ corresponds to $(\Omega_n)_{\p}$ via local class field theory. Since both $\chi$ and $\alpha_n$ are unitary, it follows that $L_{\p}(1/2, \pi_E\otimes \Omega_n) \asymp 1$ uniformly on $n$.

Finally, recall that the adjoint representation of $\tau:W_{\Q_2} \rightarrow \GL_2(\C)$ can be realized as follows. Let $\{e_1, e_2\}$ be the standard basis for $\C^2$ and let $\{e_{1}^{*}, e_{2}^{*}\}$ be its dual basis in $(\C^2)^{*}$. Then the adjoint representation of $\GL_2(\C)$ can be realized as the $3$-dimensional invariant subspace $V$ of $\C^2 \otimes (\C^2)^{*}$ that is orthogonal to the subspace generated by $e_1\otimes e_1^{*} + e_2 \otimes e_2^{*}$. (This $1$-dimensional subspace is a copy of the trivial representation.) Thus, a basis for $V$ is given by $$ \{e_1\otimes e_1^{*} - e_2\otimes e_2^{*}, e_1\otimes e_2^{*}, e_2\otimes e_1^{*} \}.$$ Under this basis, we have that $$ \mathrm{Ad}(\tau)(\Phi) = \begin{pmatrix}
    1 & 0 & 0 \\ 0 & 2^{-1} & 0 \\ 0 & 0 & 2
\end{pmatrix} $$ and, by \cite[\S 31.2]{localLanglandsGL2}, we have that the nilpotent endomorphism $N'$ corresponding to $\mathrm{Ad}(\tau)$ is given by the restriction of $N\otimes \mathrm{Id}_{2}^{*} - \mathrm{Id}_{2}\otimes N^{t}$ to $V$, where $N^{t}:(\C^{2})^{*} \rightarrow (\C^2)^{*}$ is given by $N^{t}(f) = f\circ N$. Under the given basis, we have $$ N' = \begin{pmatrix}
    0 & 0 & 1 \\ -2 & 0 & 0 \\ 0 & 0 & 0
\end{pmatrix}.$$ It follows that $\mathrm{Ker}(N')$ is the subspace of $V$ generated by $e_1\otimes e_2^{*}$, in which $\mathrm{Ad}(\tau)(\Phi)$ acts by multiplication by $2^{-1}$. Therefore, $$ L_2(1,\mathrm{Ad}(\pi)) = \mathrm{det}\left( 1-\mathrm{Ad}(\tau)(\Phi)|_{\mathrm{Ker}(N')} 2^{-1} \right)^{-1} = (1-2^{-2})^{-1} = 4/3.$$

By Proposition \ref{TwistAutoRepEtoF} we have $$ L\left(\frac{1}{2}, \pi_E \otimes \Omega_n \right) =  L\left(\frac{1}{2}, \pi \otimes \mathcal{AI}_{E/\Q}(\Omega_n)\right) = L\left( \frac{1}{2}, f_{\psi}\otimes f_{\Omega_n} \right),$$ where the $L$-function on the right hand side is the Rankin--Selberg $L$-function over $\Q$ of $f_{\psi}$ and the CM form $f_{\Omega_n}$ associated to $\Omega_n$ in Definition \ref{Def-dihedralMaass}. From the previous analysis and \eqref{PeriodIntToLvalues}, it follows that $$ I(\psi, \lambda, \rho, n) \ll_{D} \frac{Y_{\ell}^{n}(0)^2}{2\ell+1} \frac{L(\frac{1}{2},f_{\psi}\otimes f_{\Omega_n})}{L(1,\mathrm{Ad}(f_{\psi}))},$$ where the implied constants are bounded as $n$ varies. For $n\in \Z$, define $$ G(n) := \begin{cases}
    Y_{\ell}^{|n|}(0)^2/(2\ell+1) & \text{ if } |n|\leq \ell, \\
    0 & \text{ otherwise}.
\end{cases}$$ Thus, by Theorem \ref{AdjointBound}, Corollary \ref{L2FormulaFinite}, and \eqref{PeriodIntToLvalues}, we obtain 
\begin{equation}
    ||\psi|_{\mathcal{C}_{D}}||_{L^2}^{2} \ll_{D} \frac{1}{2} \sum_{\lambda \in \widehat{\Z/c_D \Z}} \sum_{\rho \in \widehat{\mathrm{Cl}(E)}} \sum_{n=0}^{\ell} \frac{Y_{\ell}^{n}(0)^2}{2\ell+1} \frac{L(\frac{1}{2},f_{\psi} \otimes f_{\Omega_n})}{L(1,\mathrm{Ad}(f_{\psi}))} \ll_{D,\varepsilon} \ell^{\varepsilon} \frac{1}{2} \sum_{n=0}^{\ell} G(n) L(1/2, f_{\psi}\otimes f_{\Omega_n}). \label{MainBound1}
\end{equation}
for any $\varepsilon>0$. We proceed to bound the latter sum.

\begin{remark}
We know that the constant $\tilde{J}_{\pi^B}(f_{\infty})$ is nonnegative for any $n\in \Z$. Since $\psi$ has real Hecke eigenvalues, it follows from \eqref{PeriodIntToLvalues} and the computations in this section that $L(\frac{1}{2},\pi_E \otimes \Omega_n) \geq 0$ for all $n\in \Z$.
\end{remark}

\subsection{Smooth partition of unity}\label{Sec-SPU}

To bound the desired sum, we start by introducing a smooth partition of unity. Let $\widetilde{U}(x)$ be a smooth bump function that is equal to $1$ on $[0,\infty)$ and supported on $[-1/2,\infty)$. Define $$ U(x) := \widetilde{U}\left( \frac{x}{\ell} \right),$$ which is equal to $1$ on $[0,\infty)$ and is supported on $[-\ell/2,\infty)$. We can see that $U(x) + U(-x) \geq 1$ for all $x\in \R$. It follows from this and $L(1/2,f_{\psi}\otimes f_{\Omega_{-n}}) = L(1/2,f_{\psi} \otimes f_{\Omega_n})$ that 
\begin{equation}
    \sum_{n=0}^{\ell} G(n) L(1/2, f_{\psi}\otimes f_{\Omega_n}) \leq \sum_{|n|\leq \ell} U(n) G(n) L(1/2,f_{\psi} \otimes f_{\Omega_n}). \label{MainBound2}
\end{equation}
Now, let $W(x)$ be the function obtained from Proposition \ref{PartitionUnity}. For each $b\geq 0$, we define the function $$ W_{b}(x) := W\left( \frac{\ell+1-x}{2^{b}} \right).$$ By construction, $W_b$ is supported on $2^{b-1} \leq \ell+1-x \leq 2^{b+1}$. We have $$ \sum_{b\geq 0} W_{b}(x) = \left\{ \begin{array}{ll}
    1 & \text{ for } x\leq \ell, \\
    0 & \text{ for } x\geq \ell+1,
\end{array} \right.$$ and is in $[0,1]$ for $\ell<x<\ell+1$. Then, from the definition of $G(n)$, we can see that
\begin{align}
    \sum_{|n|\leq \ell} U(n) G(n) L(1/2,f_{\psi} \otimes f_{\Omega_n}) &= \sum_{n\in \Z} \left(\sum_{b\geq 0} W_{b}(n) \right)U(n) G(n) L(1/2,f_{\psi} \otimes f_{\Omega_n}) \nonumber \\
    &= \sum_{b\geq 0} \sum_{n\in \Z} W_{b}(n) U(n) G(n) L(1/2,f_{\psi}\otimes f_{\Omega_n}). \label{IntroduceWell}
\end{align}
Let $\ell' = \log_2\left(\frac{3\ell}{2}+1\right)+1$. By construction, we have that $\ell+1 - 2^{b-1} \leq -\frac{\ell}{2}$ for $b \geq \ell'$, which implies that $W_b(n) U(n)=0$ for all $n\in \Z$. It follows that all terms in \eqref{IntroduceWell} with $b\geq \ell'$ vanish. Then, from this analysis and inequalities \eqref{MainBound1}, \eqref{MainBound2}, and \eqref{IntroduceWell}, it follows that
\begin{equation}
    ||\psi|_{\mathcal{C}_{D}}||_{L^2}^{2} \ll_{\varepsilon} \ell^{2\varepsilon} \sup_{0\leq b \leq \ell'} \sum_{n\in \Z} W_{b}(n) U(n) G(n) L(1/2,f_{\psi}\otimes f_{\Omega_n}). \label{MainBound3}
\end{equation}
For each $0\leq b \leq \ell'$, the inner sum on the right hand side of \eqref{MainBound3} has only finitely many terms (by the definition of $G(n)$). For the rest of the terms, the analytic conductor of the $L$-functions involved have similar analytic conductors in size, which will be useful for us later. In more rigorous terms, we have the following.

\begin{lemma}\label{AnalyticConductorIndn}
Let $0\leq b \leq \ell'$ and define $T:=2^{b}$. Then 
\begin{equation}
    G(n) \ll (T\ell)^{-1/2} \label{AnalyticConductorIndn-eq0.1}
\end{equation}
and \begin{equation}
    \mathfrak{q}_{\infty}(1/2,f_{\psi} \otimes f_{\Omega_n}) \asymp_{D} (T\ell)^2 \label{AnalyticCondunctorIndn-eq0.2}
\end{equation}
for any $n\geq -\ell/2$ such that $T/2 \leq \ell+1-n \leq 2T$.
\end{lemma}

\begin{proof}
Assume $n\geq -\ell/2$ such that $T/2 \leq \ell+1-n \leq 2T$. Since $f_{\psi}$ has level $2$ and $f_{\Omega_n}$ has level $D$, we know that $q(f_{\psi}\otimes f_{\Omega_n}) \asymp_{D} 1$, so that 
\begin{equation}
    \mathfrak{q}_{\infty}(1/2,f_{\psi}\otimes f_{\Omega_n}) \asymp_{D} (4+\ell+n)(4+\ell-n)(5+\ell+n)(5+\ell-n). \label{AnalyticConductorIndn-eq1}
\end{equation}
Now, from the inequalities we assumed, we have that $\ell-n \asymp T$. Moreover, we have that $\ell+1-2T \leq n \leq \ell+1-T/2$. Since $T \ll \ell$, we have that $n\asymp \ell$, so that $\ell+n \asymp \ell$. From these observations, \eqref{AnalyticConductorIndn-eq0.1} follows from Corollary \ref{EstimateArchimFactor}, whereas \eqref{AnalyticCondunctorIndn-eq0.2} follows from \eqref{AnalyticConductorIndn-eq1}.
\end{proof}

\begin{remark}
Indeed, the integers $n$ that satisfy $n\geq -\ell/2$ and $T/2 \leq \ell+1-n \leq 2T$ are the only ones of interest, since the support of $U(x)$ is $[-\ell/2,\infty)$ and the support of $W_{b}(x)$ is $T/2 \leq \ell+1-x \leq 2T$.
\end{remark}

Let $0\leq b \leq \ell'$ and denote $T=2^{b}$. From Lemma \ref{AnalyticConductorIndn}, it follows that 
\begin{equation}
    \sum_{n\in \Z} W_{b}(n) U(n) G(n) L(1/2,f_{\psi}\otimes f_{\Omega_n}) \ll (T\ell)^{-1/2} \sum_{n\in \Z} W_{b}(n) U(n) L(1/2,f_{\psi}\otimes f_{\Omega_n}). \label{MainBound4.5}
\end{equation}
We now proceed to prove that, for any $\varepsilon>0$ and $0\leq b \leq \ell'$, we have
\begin{equation*}
    (T\ell)^{-1/2} \sum_{n\in \Z} W_{b}(n) U(n) L(1/2,f_{\psi}\otimes f_{\Omega_n}) \ll_{D,\varepsilon} \ell^{\varepsilon},
\end{equation*}
which would imply the main result.

\subsection{Approximate functional equation}

Before using the approximate functional equation in our case, we introduce some notation. 

\begin{definition}\label{DefVs}
For $n\in \Z$, we let $V_s(y,n)$ be the function defined in Theorem \ref{AppFuncEq1} for $f=f_{\psi}$ and $g=f_{\Omega_n}$. More explicitly, $$ V_s(y,n) := \frac{1}{2\pi i} \int_{\mathrm{Re}(u)=3} \frac{\gamma(s+u,f_{\psi}\otimes f_{\Omega_n})}{\gamma(s,f_{\psi}\otimes f_{\Omega_n})} y^{-u} e^{u^2} \: \frac{\mathrm{d}u}{u}.$$
\end{definition}

\begin{remark}\label{VsynSmoothWRTn}
In fact, since $f_{\psi}$ has weight $2\ell+2$ and $f_{\Omega_n}$ has weight $2n+1$, we have that 
\begin{align*}
    &\gamma(s,f_{\psi}\otimes f_{\Omega_n}) \\
    &= \pi^{-2s} \Gamma\left( \frac{s+\ell-n+\frac{1}{2}}{2} \right) \Gamma\left( \frac{s+\ell-n+\frac{3}{2}}{2} \right) \Gamma\left( \frac{s+\ell+n+\frac{1}{2}}{2} \right) \Gamma\left( \frac{s+\ell+n+\frac{3}{2}}{2} \right).
\end{align*}
It follows that the functions $\gamma(s,f_{\psi} \otimes f_{\Omega_n})$ and $V_{s}(y,n)$ may be extended to smooth functions $\gamma(s,f_{\psi},x)$ and $V_{s}(y,x)$ of $x$ on their domain. We make use of this fact later.
\end{remark}

Note that $\bigo_{E}^{\times}$ acts on $\bigo_{E}$ by left multiplication, and that two elements of $\bigo_{E}$ generate the same ideal in $\bigo_{E}$ if and only if they are in the same orbit of $\bigo_{E}^{\times} \backslash \bigo_{E}$.

\begin{definition}
Let $\mathcal{F}_{D}$ be defined by $$ \mathcal{F}_{D} := \{0\} \cup \left\{ \alpha \in \bigo_{E} : \mathrm{Re}(\alpha) > 0 \right\} \cup \left\{ \alpha \in \bigo_{E} : \mathrm{Re}(\alpha)=0 \text{ and } \mathrm{Im}(\alpha) <0 \right\},$$ which is a fundamental domain for the action of $\bigo_{E}^{\times}$ on $\bigo_{E}$.
\end{definition}

\begin{remark}
The fundamental domain $\mathcal{F}_{D}$ is different for $D=-4$ and $D=-3$. For instance, for $D=-4$ we have $$ \mathcal{F}_{-4} = \{0\} \cup \left\{ \alpha \in \bigo_{E} : \mathrm{Im}(\alpha) \geq 0 \text{ and } \mathrm{Re}(\alpha) >0 \right\}.$$
\end{remark}

To be able to use the approximate functional equation for the central values of the $L$-functions involved, we need an estimate of the analytic conductors, which was done in Lemma \ref{AnalyticConductorIndn}.

\begin{proposition}
For any $b$ fixed, $\delta>0$ and $B\geq 1$, we have
\begin{align}
    &(T\ell)^{-1/2}\sum_{n\in \Z} W_{b}(n) U(n) L(1/2,f_{\psi}\otimes f_{\Omega_n}) \label{MainBound4} \\
    &= 2\sum_{m=1}^{(T\ell)^{1+\delta}} \frac{\chi_{D}(m)}{m} \sum_{\substack{\alpha \in \mathcal{F}_{D} \\ 1\leq N(\alpha) \leq (T\ell)^{1+\delta}}} \frac{\lambda_{f_{\psi}}(N(\alpha))}{\sqrt{N(\alpha)}} \Pi_{\alpha}(m) + O_{D,B,\delta}\left(\ell^{-2B}\right), \nonumber
\end{align}
where $$ \Pi_{\alpha}(m) := (T\ell)^{-1/2}\sum_{n\in \Z} \Omega_n((m)) W_{b}(n) U(n) V_{1/2}\left(\frac{m^2 N(\alpha)}{\sqrt{q(f_{\psi}\otimes f_{\Omega_n})}},n\right) \Omega_n((\alpha))$$ and $V_{s}(y,n)$ is the function defined in Definition \ref{DefVs}.
\end{proposition}

\begin{proof}
It follows from Lemma \ref{AnalyticConductorIndn} and the approximate functional equation (Remark \ref{AppFuncEqHoldsForDeven} and Corollary \ref{AppFuncEq2}) for $\delta>0$, $M=(T\ell)^{1+\delta}$, and $c$ the inverse of the implicit constant in the upper bound of Lemma \ref{AnalyticConductorIndn}.
\end{proof}

We want to estimate the main term on the right hand side of \eqref{MainBound4}. For this, we again introduce a smooth dyadic partition of unity. Let $a_T$ be the smallest integer with $a_T > (T\ell)^{1+\delta}$ and let $U_T$ be a smooth bump function supported on $[0,a_T]$ and equal to $1$ on $[1,(T\ell)^{1+\delta}]$. For $\alpha \in \bigo_K$, by abuse of notation we write $$ U_T(\alpha) := U_T(N(\alpha)).$$ Let $W(x)$ be the function from Proposition \ref{PartitionUnity} and, for $a\geq 0$ and $\alpha \in \bigo_K$, define $$ W_{a}' := W\left( \frac{N(\alpha)}{2^{a}} \right).$$ 
Then, similar to Section \ref{Sec-SPU}, for all $1\leq m\leq (T\ell)^{1+\delta}$ we have 
\begin{equation}
    \sum_{\substack{\alpha \in \mathcal{F}_{D} \\ 1\leq N(\alpha) \leq (T\ell)^{1+\delta}}} \frac{\lambda_{f_{\psi}}(N(\alpha))}{\sqrt{N(\alpha)}} \Pi_{\alpha}(m) = \sum_{a\geq 0} \sum_{\alpha \in \mathcal{F}_{D}} \frac{\lambda_{f_{\psi}}(N(\alpha))}{\sqrt{N(\alpha)}} \Pi_{\alpha}(m) W_{a}'(\alpha) U_T(\alpha). \label{SecondPartitionUnity}
\end{equation}
Again, as in Section \ref{Sec-SPU}, whenever $2^{a-1} \geq a_T$, we have $W_{a}'(\alpha) U_T(\alpha)=0$ for all $\alpha \in \mathcal{F}_{D}$, so the $a$-sum on the right hand side of \eqref{SecondPartitionUnity} consists of $\ll_{\delta} \log(T\ell)$ terms. Since $T\ll \ell$, from \eqref{MainBound4} and \eqref{SecondPartitionUnity} it follows that
\begin{align}
    &(T\ell)^{-1/2} \sum_{n\in \Z} W_{b}(n) U(n) L(1/2,f_{\psi}\otimes f_{\Omega_n}) \label{MainBound5} \\
    &\ll_{D, \varepsilon, \delta} \ell^{\varepsilon} \sup_{0\leq 2^{a} \ll (T\ell)^{1+\delta}} \left| \sum_{m=1}^{(T\ell)^{1+\delta}} \frac{\chi_{D}(m)}{m} \sum_{\alpha \in \mathcal{F}_{D}} \frac{\lambda_{f_{\psi}}(N(\alpha))}{\sqrt{N(\alpha)}} \Pi_{\alpha}(m) W_{a}'(\alpha) U_T(\alpha) \right| \nonumber
\end{align}
for all $\varepsilon>0$, $\delta>0$, and $1\leq T \ll \ell$.

\subsection{Partial derivatives of \texorpdfstring{$V_{1/2}$}{V1/2}}

To proceed, we will need certain estimates on the partial derivatives of $V_{1/2}(y,x)$ with respect to $x$. We show these estimates in this subsection. First we estimate the quotient of gamma factors inside the integral in the definition of $V_{1/2}(y,x)$ and its derivatives.

\begin{lemma}\label{BoundsGammas}
Let $\delta >0$ and let $1\leq T \ll \ell$. Suppose that $T/2 \leq \ell+1-x \leq 2T$ and $x\geq -\ell/2$. Define $F(u,f_{\psi},x) := \frac{\gamma(u+1/2,f_{\psi},x)}{\gamma(1/2,f_{\psi},x)}$, where $$ \gamma(s,f_{\psi},x) := \pi^{-2s} \Gamma\left( \frac{s+\ell-x+\frac{1}{2}}{2} \right) \Gamma\left( \frac{s+\ell-x+\frac{3}{2}}{2} \right) \Gamma\left( \frac{s+\ell+x+\frac{1}{2}}{2} \right) \Gamma\left( \frac{s+\ell+x+\frac{3}{2}}{2} \right).$$ Then, for $\mathrm{Re}(u) = \delta$, we have $$ |F(u,f_{\psi},x)| \ll_{\delta} e^{3\pi |\mathrm{Im}(u)|} \ell^{2\delta}.$$
\end{lemma}

\begin{proof}
Recall Stirling's estimate (see \cite[Thm. C.1]{MontgomeryVaughan})
\begin{equation}
    |\Gamma(s)| \asymp \left|s^{s-1/2}\right| e^{-\mathrm{Re}(s)}, \label{StirlingGamma}
\end{equation}
which holds for $|s|\geq \delta$ and $|\mathrm{Arg}(s)| < \pi-\delta$. Note that $$ \left| s^{s-1/2} \right| = |s|^{\mathrm{Re}(s)-1/2} e^{-\mathrm{Arg}(s) \mathrm{Im}(s)}.$$ Using this observation and \eqref{StirlingGamma}, we have
\begin{equation}
    |F(u,f_{\psi},x)| \asymp F_1 F_2 F_3, \label{BoundsGammas-1}
\end{equation}
where
\begin{align*}
    F_1 &= \frac{|\ell-x+1+u|^{(\ell-x)/2} |\ell-x+2+u|^{(\ell-x+1)/2} }{(\ell-x+1)^{(\ell-x)/2} (\ell-x+2)^{(\ell-x+1)/2} } \\
    &\quad \times \frac{|\ell+x+1+u|^{(\ell+x)/2} |\ell+x+2+u|^{(\ell+x+1)/2}}{(\ell+x+1)^{(\ell+x)/2} (\ell+x+2)^{(\ell+x+1)/2}}, \\
    F_2 &= |\ell-x+1+u|^{\delta/2} |\ell-x+2+u|^{\delta/2} |\ell+x+1+u|^{\delta/2} |\ell+x+2+u|^{\delta/2}, \\
    F_3 &= e^{-\frac{1}{2}\mathrm{Im}(u) \left[ \mathrm{Arg}((\ell-x+1+u)/2) + \mathrm{Arg}((\ell-x+2+u)/2) + \mathrm{Arg}((\ell+x+1+u)/2) + \mathrm{Arg}((\ell+x+2+u)/2) \right]}.
\end{align*}
We proceed to estimate each factor $F_i$. We start with $F_1$. By the triangle inequality, we have
\begin{equation}
    F_{1}^2 \leq \left( 1+\frac{|u|}{\ell-x+1} \right)^{\ell-x} \left( 1+\frac{|u|}{\ell-x+2} \right)^{\ell-x+1} \left( 1+\frac{|u|}{\ell+x+1} \right)^{\ell+x} \left( 1+\frac{|u|}{\ell+x+2} \right)^{\ell+x+1}. \label{BoundsGammas-2}
\end{equation}
For any $C\geq 0$, the function $f(t) = \left( 1+\frac{C}{t+1}\right)^t$ is nondecreasing on the interval $[-1,\infty)$ and is bounded from above by $e^C$. From this observation and \eqref{BoundsGammas-2}, we obtain
\begin{equation}
    F_{1} \leq e^{2|u|} \ll_{\delta} e^{\pi |\mathrm{Im}(u)|}. \label{BoundsGammas-3}
\end{equation}

We now estimate $F_2$. First, by definition of $x$ and $T$, we have that $\ell-x+1\ll \ell$ and $\ell+x+1 \ll \ell$. Since $\ell\geq 1$, we have that $ |\ell\pm x+1+u| \ll \ell+|u| = \ell(1+|u|/\ell) \leq \ell(1+|u|)$. It follows, since $\mathrm{Re}(u)=\delta$, that 
\begin{equation}
    F_2 \ll \ell^{2\delta} (1+|u|)^{2\delta} \ll_{\delta} \ell^{2\delta} e^{\pi |\mathrm{Im}(u)|}. \label{BoundsGammas-4}
\end{equation}

Finally, we estimate $F_3$. Notice that, by definition of $x$ and $T$, we have that $\ell-x+1>0$ and $\ell+x+1 >0$. Hence, each argument in the exponential defining $F_3$ is between $-\pi/2$ and $\pi/2$, so in absolute value they are $\leq \pi/2$. Thus
\begin{equation}
    F_3 \leq e^{\frac{1}{2}|\mathrm{Im}(u)| \left[ \pi/2 + \pi/2 + \pi/2 + \pi/2\right]} = e^{\pi|\mathrm{Im}(u)|}. \label{BoundsGammas-5}
\end{equation}
Combining \eqref{BoundsGammas-1}, \eqref{BoundsGammas-3}, \eqref{BoundsGammas-4}, and \eqref{BoundsGammas-5} the result follows.
\end{proof}

To estimate the derivatives of $F$, we need some preliminary results. These are in \cite[\S 3.2.1]{DanaMScThesis} and we reproduce them here for the ease of exposition. The first one of them is a consequence of the mean value theorem.

\begin{lemma}\label{MVTComplex}
Let $\Omega \subseteq \C$ be an open subset containing $0$ and let $f:\Omega \rightarrow \C$ be a holomorphic function. Let $z_0 \in \Omega$ be such that the line segment joining $0$ and $z_0$ is contained in $\Omega$. Then $$ |f(z_0) - f(0)| \leq 2|z_0| \sup_{r\in [0,1]} |f'(rz_0)|.$$
\end{lemma}

Recall that, for $m\geq 0$ an integer, the polygamma function $\psi^{(m)}$ is the meromorphic function defined as $$ \psi^{(m)}(z) := \left\{ \begin{array}{ll}
    \Gamma'(z)/\Gamma(z) & \text{ if } m=0, \\
    \frac{\mathrm{d}^m}{\mathrm{d}z^m} \psi^{(0)}(z) & \text{ if } m\geq 1.
\end{array} \right.$$ Using the series representation for $\psi^{(m)}$ shown in \cite[6.4.10]{TableFormulas}, we have the following estimate.

\begin{proposition}\label{PolygammaEstimates}
For $\mathrm{Re}(z)>0$ and $m\geq 1$, we have $$ \psi^{(m)}(z) \ll_{m} |z|^{-m-1} + |z|^{-m}.$$
\end{proposition}

The polygamma functions allow us to write the derivatives of the gamma function as polynomials in these functions. The same is true for quotients of gamma functions, which is what we have in our case.

\begin{definition}
For $m\geq 0$, $s\in \C$, and $x\in \R$, we let $\psi^{(m)}(s,f_{\psi},x)$ be the expression $$\psi^{(m)}\left( \frac{s+\ell+x+\frac{1}{2}}{2} \right) +\psi^{(m)}\left( \frac{s+\ell+x+\frac{3}{2}}{2} \right) -\psi^{(m)}\left( \frac{s+\ell-x+\frac{1}{2}}{2} \right) -\psi^{(m)}\left( \frac{s+\ell-x+\frac{3}{2}}{2} \right).$$
\end{definition}

\begin{proposition}\label{Derivativegamma}
For $x\in \R$ and $s\in \C$, we have $$ \frac{\partial}{\partial x} \gamma(s,f_{\psi},x) = \frac{1}{2} \psi^{(0)}(s,f_{\psi},x) \gamma(s,f_{\psi},x).$$
\end{proposition}

Moreover, we have an analogue of Proposition \ref{PolygammaEstimates} for $\psi^{(m)}(s,f_{\psi},x)$.

\begin{proposition}\label{ModifiedPolygammaEstimates}
Let $1\leq T \ll \ell$. Assume that $T/2 \leq \ell+1-x \leq 2T$ and $x\geq -\ell/2$. Then, for any $m\geq 1$ and $\mathrm{Re}(s)\geq 1/2$, we have $$ \psi^{(m)}(s,f_{\psi},x) \ll_{m} T^{-m}.$$
\end{proposition}

With the previous results, we can bound the partial derivatives of $F$ with respect to $x$.

\begin{proposition}\label{DerivativesF}
Let $1\leq T \ll \ell$ and $\delta >0$. Suppose that $T/2 \leq \ell+1-x \leq 2T$,  $x\geq -\ell/2$, and $\mathrm{Re}(u)=\delta$. Then, for any $j\geq 0$ we have $$ \frac{\partial^j}{\partial x^j} F(u,f_{\psi},x) \ll_{\delta, j} |u|^{j} e^{3\pi |\mathrm{Im}(u)|} \ell^{2\delta} T^{-j}.$$
\end{proposition}

\begin{proof}
We will show that
\begin{equation}
    \frac{\partial^j}{\partial x^j} F(u,f_{\psi},x) \ll_{\delta, j} |u|^{j} \: |F(u,f_{\psi},x)|. \label{DerivativesF-1}
\end{equation}
Then the result will follow from Lemma \ref{BoundsGammas}. First, define $$ D(u,f_{\psi},x) := \psi^{(0)}(1/2 + u,f_{\psi},x) - \psi^{(0)}(1/2,f_{\psi},x).$$ By Proposition \ref{Derivativegamma}, we can see that 
\begin{align*}
    \frac{\partial}{\partial x} F(u,f_{\psi},x) &= \frac{\gamma(1/2,f_{\psi},x) \cdot \frac{\partial}{\partial x} \gamma(1/2+u,f_{\psi},x) - \gamma(1/2+u,f_{\psi},x) \cdot \frac{\partial}{\partial x} \gamma(1/2,f_{\psi},x)}{\gamma(1/2, f_{\psi},x)^2} \\
    &= \frac{1}{2} D(u,f_{\psi},x) F(u,f_{\psi},x).
\end{align*}
We first estimate the partial derivatives of $D$ with respect to $x$. First, note that for any $j\geq 0$, we have $$ \frac{\partial^{j}}{\partial x^{j}} D(u,f_{\psi},x) = \frac{1}{2^j}\left( \psi^{(j)}(1/2+u,f_{\psi},x) - \psi^{(j)}(1/2,f_{\psi},x) \right).$$ Hence, by Lemma \ref{MVTComplex} and Proposition \ref{ModifiedPolygammaEstimates}, we have that
\begin{equation*}
    \frac{\partial^{j}}{\partial x^{j}} D(u,f_{\psi},x) \ll_{j} |u| \: \sup_{r\in [0,1]} \left|\psi^{(j+1)}(1/2 + ru, f_{\psi},x)\right| \ll_{j} |u| T^{-j-1}. 
\end{equation*}
We now come back to estimating the partial derivatives of $F$ with respect to $x$. To show \eqref{DerivativesF-1}, we proceed by induction on $j$. Notice that \eqref{DerivativesF-1} holds for $j=0$ and $j=1$. Now, assume that \eqref{DerivativesF-1} holds for all $m\leq j$. By the induction hypothesis, we have that 
\begin{align*}
    \frac{\partial^{j+1}}{\partial x^{j+1}} F(u,f_{\psi},x) &\ll_{j} \sum_{m=0}^{j} \left| \frac{\partial^{j-m}}{\partial x^{j-m}} F(u,f_{\psi},x) \right| \left| \frac{\partial^{m}}{\partial x^{m}} D(u,f_{\psi},x) \right| \\
    &\ll_{j} \sum_{m=0}^{j} |u|^{j-m} \: |F(u,f_{\psi},x)| T^{-j+m} \cdot |u|T^{-m-1} \\
    &\ll_{j} |u|^{j+1} \: |F(u,f_{\psi},x)| T^{-j-1},
\end{align*}
which completes the proof.
\end{proof}

From these bounds, we can estimate the partial derivatives of $V_{1/2}(y,x)$ with respect to $x$.

\begin{proposition}\label{BoundsPartialsVs}
Let $1\leq T \ll \ell$. Suppose that $T/2 \leq \ell+1-x \leq 2T$ and $x\geq -\ell/2$. Then, for any $j\geq 1$ and $\varepsilon>0$, we have $$ \frac{\partial^{j}}{\partial x^{j}} V_{1/2}(y,x) \ll_{\varepsilon,j} y^{-\varepsilon/2} \ell^{\varepsilon} T^{-j}.$$    
\end{proposition}
    
\begin{proof}
Recall from Definition \ref{DefVs} that 
\begin{equation}
    V_{1/2}(y,x) := \frac{1}{2\pi i} \int_{\mathrm{Re}(u)=3} F(u,f_{\psi},x) y^{-u} e^{u^2} \: \frac{\mathrm{d}u}{u} \label{BoundsPartialsVs-1}
\end{equation}
Let $\delta = \varepsilon/2$. Since $e^{u^2}$ decays like $e^{-\mathrm{Im}(u)^2}$ on vertical lines, by Stirling's estimate of the Gamma function and Cauchy's theorem, we can move the line of integration in \eqref{BoundsPartialsVs-1} to the line $\mathrm{Re}(u) = \delta$, so that 
\begin{equation*}
    V_{1/2}(y,x) = \frac{1}{2\pi i} \int_{\mathrm{Re}(u)=\delta} F(u,f_{\psi},x) y^{-u} e^{u^2} \: \frac{\mathrm{d}u}{u}.
\end{equation*}
It follows by Proposition \ref{DerivativesF} that
\begin{align*}
    \frac{\partial^{j}}{\partial x^{j}} V_{1/2}(y,x) &\ll \int_{\mathrm{Re}(u)=\delta} \left| \frac{\partial^{j}}{\partial x^{j}} F(u,f_{\psi},x) \right| y^{-\mathrm{Re}(u)} e^{\mathrm{Re}(u^2)} \frac{\mathrm{d}u}{|u|} \\
    &\ll_{\delta,j} y^{-\delta} e^{\delta^2} \ell^{2\delta} T^{-j} \int_{\mathrm{Re}(u)=\delta} |u|^{j-1} e^{3\pi |\mathrm{Im}(u)|} e^{-|\mathrm{Im}(u)|^2} \: \mathrm{d}u \\
    &\ll_{\delta,j} y^{-\delta} \ell^{2\delta} T^{-j}. \qedhere
\end{align*}
\end{proof}

\subsection{Poisson summation}

Notice that, by the construction of $\Omega_{n}$, for $m\in \Z$ and $\alpha \in \bigo_E$, we have that $$ \Omega_{n}((m)) = \left(\frac{m}{\overline{m}}\right)^{n} = 1 \quad \text{and} \quad \Omega_{n}((\alpha)) = \left(\frac{\alpha}{\overline{\alpha}}\right)^{n} = e^{in\cdot \mathrm{Arg}(\alpha/\overline{\alpha})}.$$ Let
$$ q_D := \left\{ \begin{array}{ll}
    2|D| & \text{ if } D \text{ is odd}, \\
    \sqrt{2}|D| & \text{ if } D \text{ is even}.
\end{array} \right. $$ By \cite[Thm. 2.2]{WinnieLiRankinSelberg}, we have that $q(f_{\psi}\otimes f_{\Omega_{n}})=q_{D}^2$ for all $n\in \Z$. Therefore, $$ \Pi_{\alpha}(m) = (T\ell)^{-1/2} \sum_{n\in \Z} W_{b}(n) U(n) V_{1/2}\left(\frac{m^2 N(\alpha)}{q_D},n\right) e^{in\cdot \mathrm{Arg}(\alpha/\overline{\alpha})}. $$

We now fix $0\leq a\leq a_T$ and write $A:= 2^{a}$. We want to estimate the expression inside the absolute value on the right hand side of \eqref{MainBound5}, which we call $S_a$. Define the smooth compactly supported function $$ H_A(\alpha,m,x) := \frac{(T\ell)^{-1/2}}{\sqrt{N(\alpha)}} W_{a}'(\alpha) U_T(\alpha) W_{b}(x) U(x) V_{1/2}\left( \frac{m^2 N(\alpha)}{q_D},x \right).$$
Thus, we have $$ S_a = \sum_{m=1}^{(T\ell)^{1+\delta}} \frac{\chi_{D}(m)}{m} \sum_{\alpha \in \mathcal{F}_{D}} \lambda_{f_{\psi}}(N(\alpha)) \sum_{n\in \Z} H_A(\alpha,m,n) \, e^{in\cdot \mathrm{Arg}(\alpha/\overline{\alpha})}.$$ 
Consequently, we can use the Poisson summation formula in the $n$-sum to obtain that $$ S_a = \sum_{m=1}^{(T\ell)^{1+\delta}} \frac{\chi_{D}(m)}{m} \sum_{\alpha \in \mathcal{F}_{D}} \lambda_{f_{\psi}}(N(\alpha)) \sum_{\xi\in \Z} \widehat{H_A}\left(\alpha,m,\xi - \frac{\mathrm{Arg}(\alpha/\overline{\alpha})}{2\pi} \right),$$ where $\widehat{H_A}$ is the Fourier transform of $H_A$ with respect to $x$. 

\begin{remark}
For $D=-3$ and $D=-4$, one has to split the $n$-sum in the definition of $S_a$ into several sums depending on the congruence of $n$ modulo $3$ or $4$, respectively. There are two reasons for this change. The first one is the difference in the definition of $\Omega_n$ for different values of $n$, as pointed out in Remark \ref{D34HeckeCharacters}. The second one is that, along each sum, we want $q(f_{\psi} \otimes f_{\Omega_n})$ to remain constant, so that we can use the Poisson summation formula. However, each sum can be bounded following the procedure we describe now.
\end{remark}

For $\xi\neq 0$, by repeated use of integration by parts, we have
\begin{equation}
    \widehat{H_A}(\alpha,m,\xi) \leq \frac{1}{(2\pi |\xi|)^{j}} \left| \widehat{H_{A}^{(j)}}(\alpha,m,\xi) \right|, \label{BoundFourierTransDer}
\end{equation}
where $H_A^{(j)}$ is the $j$-th partial derivative of $H_A(\alpha,m,x)$ with respect to $x$. We estimate $\widehat{H_A}$ via estimating $H_{A}^{(j)}(\alpha,m,x)$.

\begin{proposition}\label{BoundDerivativesHA}
Let $j\geq 0$. Then $$ H_{A}^{(j)}(\alpha,m,x) = 0 \text{  unless  } 2x\in \left[ \ell+1-2T, \ell+1-T/2 \right]\cap [-\ell/2,\infty)  \text{  and  } N(\alpha) \in [A/2,2A].$$ Moreover, for any $\varepsilon>0$, we have $$ H_{A}^{(j)}(\alpha,m,x) \ll_{D,\varepsilon,j} \ell^{\varepsilon} T^{-j} (\ell TA)^{-1/2}.$$
\end{proposition}

\noindent \textit{Proof.} The first statement follows from the definitions of $W_{a}'(\alpha)$ and $W_{b}(x)$, namely that $W_{a}'$ has support on $N(\alpha) \in [A/2,2A]$ and that $W_{b}$ has support on $[\ell+1-2T,\ell+1-T/2]$. Now, assume that $N(\alpha) \in [A/2,2A]$, so that $N(\alpha) \gg A$. By the product rule for derivatives, we have
\begin{align}
    H_{A}^{(j)}(\alpha,m,x) &\ll (\ell TA)^{-1/2} \left( W_{b}(x) U(x) V_{1/2}\left( \frac{m^2 N(\alpha)}{q_D},x \right) \right)^{(j)} \label{BoundDerivativesHA-1} \\
    &\ll_{j} (\ell TA)^{-1/2} \sum_{i=0}^{j} \left(W_{b}(x) U(x)\right)^{(j-i)} V_{1/2}^{(i)}\left( \frac{n^2 N(\alpha)}{q_D},x \right), \nonumber
\end{align}
where all the derivatives are with respect to $x$. Since $W_{b}(x) = W\left(\frac{\ell-x}{T}\right)$ for $W(x)$ a fixed smooth compactly supported function, by using the chain rule several times we obtain that $W_{b}^{(j-i)}(2x) \ll_{j,i} T^{-j+i}$. Similarly, since $U(x) = \widetilde{U}\left(\frac{x}{\ell}\right)$ for $\widetilde{U}(x)$ a fixed bump function, we have that $U^{(j-i)}(2x) \ll_{j,i} \ell^{-j+i} \ll T^{-j+i}$. It follows that 
\begin{equation}
    \left(W_{b}(2x) U(2x) \right)^{(j-i)} \ll_{j,i} \sum_{h=0}^{j-i} W_{b}^{(h)}(2x) U^{(j-i-h)}(2x) \ll_{j,i} \sum_{h=0}^{j-i} T^{-h} T^{h+i-j} \ll_{j,i} T^{-j+i}. \label{BoundDerivativesHA-2}
\end{equation}
Additionally, from Proposition \ref{BoundsPartialsVs}, we have that
\begin{equation}
    V_{1/2}^{(i)}\left(\frac{m^2 N(\alpha)}{q_D},x\right) \ll_{\varepsilon,i} q_{D}^{\varepsilon/2} \ell^{\varepsilon} T^{-i}. \label{BoundDerivativesHA-3}
\end{equation}
From \eqref{BoundDerivativesHA-1}, \eqref{BoundDerivativesHA-2}, and \eqref{BoundDerivativesHA-3}, we obtain
\[
\pushQED{\qed}
H_{A}^{(j)}(\alpha,m,x) \ll_{D,\varepsilon,j} (\ell TA)^{-1/2} \sum_{i=0}^{j} T^{-j+i} \ell^{\varepsilon} T^{-i} \ll_{D,\varepsilon,j} (\ell TA)^{-1/2} \ell^{\varepsilon} T^{-j}. \qedhere
\popQED
\]

As a consequence of \eqref{BoundFourierTransDer} and Proposition \ref{BoundDerivativesHA}, we have the following.

\begin{corollary}\label{BoundHATransform}
For any $j\geq 0$ and $\xi\neq 0$, we have $$ \widehat{H_A}(\alpha,m,\xi) \ll_{D,\varepsilon,j} \frac{T \ell^{\varepsilon}}{|T\xi|^{j}} (\ell TA)^{-1/2}.$$ Moreover, for $N(\alpha) \notin [A/2,2A]$, we have $\widehat{H_A}(\alpha,m,\xi)=0$.
\end{corollary}

Using this last estimate, we can bound the inner two sums of $S_a$.

\begin{proposition}\label{BoundInnerSumsSa}
For any $m\geq 1$ and $0 < \delta \leq 1$, we have 
\begin{align*}
    &\sum_{\alpha \in \mathcal{F}_{D}} \lambda_{f_{\psi}}(N(\alpha)) \sum_{\xi\in \Z} \widehat{H_A}\left( \alpha,m,\xi-\frac{\mathrm{Arg}(\alpha/\overline{\alpha})}{2\pi} \right) \\
    &\ll_{D,\varepsilon,\delta} \ell^{\varepsilon} T^{1/2} (\ell A)^{-1/2} \sum_{\alpha \in S_{T^{-1+\delta}}} |\lambda_{f_{\psi}}(N(\alpha))| + \ell^{\varepsilon} T^{-3/2} (\ell A)^{-1/2} \sum_{\substack{\alpha \in \mathcal{F}_{D} \\ A/2 \leq N(\alpha) \leq 2A}} |\lambda_{f_{\psi}}(N(\alpha))|,
\end{align*}
where $$ S_{R} := \left\{ \alpha \in \mathcal{F}_{D} : -\pi R \leq \mathrm{Arg}(\alpha) \leq \pi R, \: A/2 \leq N(\alpha) \leq 2A \right\}$$ for any $R\in (0,1)$.
\end{proposition}

\begin{proof}
From the definition of $\mathcal{F}_{D}$, for any $\alpha \in \mathcal{F}_{D}$ we have that $-\pi/2 \leq \mathrm{Arg}(\alpha) < \pi/2$, so that $-\pi \leq \mathrm{Arg}(\alpha/\overline{\alpha}) < \pi$. This implies that $-\frac{1}{2} \leq \frac{\mathrm{Arg}(\alpha / \overline{\alpha})}{2\pi} < \frac{1}{2}$ for any $\alpha \in \mathcal{F}_{D}$. Thus, from Corollary \ref{BoundHATransform} with $j=2$ and the triangle inequality, we have that
\begin{align}
    &\sum_{\alpha \in \mathcal{F}_{D}} \lambda_{f_{\psi}}(N(\alpha)) \sum_{\xi\geq 1} \widehat{H_A}\left( \alpha,m,\xi-\frac{\mathrm{Arg}(\alpha/\overline{\alpha})}{2 \pi} \right) \label{BoundInnerSumsSa-1} \\
    &\ll_{D,\varepsilon} \ell^{\varepsilon} T^{-3/2} (\ell A)^{-1/2} \sum_{\substack{\alpha \in \mathcal{F}_{D} \\ A/2 \leq N(\alpha) \leq 2A}} |\lambda_{f_{\psi}}(N(\alpha))| \sum_{\xi\geq 1} \frac{1}{\left| \xi - \frac{\mathrm{Arg}(\alpha/\overline{\alpha})}{2 \pi} \right|^2} \nonumber \\
    &\ll_{D,\varepsilon} \ell^{\varepsilon} T^{-3/2} (\ell A)^{-1/2} \sum_{\substack{\alpha \in \mathcal{F}_{D} \\ A/2 \leq N(\alpha) \leq 2A}} |\lambda_{f_{\psi}}(N(\alpha))| \sum_{\xi\geq 1} \frac{1}{\xi^2} \nonumber \\
    &\ll_{D,\varepsilon} \ell^{\varepsilon} T^{-3/2} (\ell A)^{-1/2} \sum_{\substack{\alpha \in \mathcal{F}_{D} \\ A/2 \leq N(\alpha) \leq 2A}} |\lambda_{f_{\psi}}(N(\alpha))|. \nonumber
\end{align}
Using the same reasoning, we have
\begin{equation}
    \sum_{\alpha \in \mathcal{F}_{D}} \lambda_{f_{\psi}}(N(\alpha)) \sum_{\xi\leq -1} \widehat{H_A}\left( \alpha,m,\xi-\frac{\mathrm{Arg}(\alpha/\overline{\alpha})}{2 \pi} \right) \ll_{D,\varepsilon} \ell^{\varepsilon} T^{-3/2} (\ell A)^{-1/2} \sum_{\substack{\alpha \in \mathcal{F}_{D} \\ A/2 \leq N(\alpha) \leq 2A}} |\lambda_{f_{\psi}}(N(\alpha))|. \label{BoundInnerSumsSa-2}
\end{equation}
Now, let $\xi = 0$ and $\alpha \in \mathcal{F}_{D}$. By the definition of $S_{R}$, we have that $\left|-\frac{\mathrm{Arg}(\alpha/\overline{\alpha})}{2 \pi} \right| \leq T^{-1+\delta}$ if and only if $\alpha \in S_{T^{-1+\delta}}$. If this is the case, from Corollary \ref{BoundHATransform} with $j=0$, we get 
\begin{equation}
    \widehat{H_A}\left(\alpha,m,-\frac{\mathrm{Arg}(\alpha/\overline{\alpha})}{2 \pi} \right) \ll_{D,\varepsilon} \ell^{\varepsilon} T^{1/2} (\ell A)^{-1/2}. \label{BoundInnerSumsSa-3}
\end{equation}
If not, then we have $\left|-\frac{\mathrm{Arg}(\alpha/\overline{\alpha})}{2 \pi} \right| > T^{-1+\delta}$. In particular, $\frac{\mathrm{Arg}(\alpha/\overline{\alpha})}{2 \pi} \neq 0$. Hence, again by Corollary \ref{BoundHATransform} with $j=\left\lceil 2/\delta \right\rceil$ we have 
\begin{align}
    \widehat{H_A}\left(\alpha,m,-\frac{\mathrm{Arg}(\alpha/\overline{\alpha})}{2\pi} \right) &\ll_{D,\varepsilon,\delta} \frac{\ell^{\varepsilon} T}{\left| T\left( -\frac{\mathrm{Arg}(\alpha/\overline{\alpha})}{\pi} \right) \right|^{j}}  (\ell TA)^{-1/2} \label{BoundInnerSumsSa-4} \\
    &\ll_{D,\varepsilon,\delta} \ell^{\varepsilon} T^{1-\delta j} (\ell TA)^{-1/2} \ll_{D,\varepsilon,\delta} \ell^{\varepsilon} T^{-3/2} (\ell A)^{-1/2}. \nonumber
\end{align}
Combining \eqref{BoundInnerSumsSa-3} and \eqref{BoundInnerSumsSa-4}, we obtain
\begin{align}
    &\sum_{\alpha \in \mathcal{F}_{D}} \lambda_{f_{\psi}}(N(\alpha)) \widehat{H_A}\left( \alpha,m,-\frac{\mathrm{Arg}(\alpha/\overline{\alpha})}{2 \pi} \right) \label{BoundInnerSumsSa-5} \\
    &\leq \sum_{\substack{\alpha \in \mathcal{F}_{D} \\ A/2 \leq N(\alpha) \leq 2A}} |\lambda_{f_{\psi}}(N(\alpha))| \cdot \left| \widehat{H_A}\left( \alpha,m,-\frac{\mathrm{Arg}(\alpha/\overline{\alpha})}{2 \pi} \right) \right| \nonumber \\
    &\ll_{D,\varepsilon, \delta} \ell^{\varepsilon} T^{1/2} (\ell A)^{-1/2} \sum_{\alpha\in S_{T^{-1+\delta}}} |\lambda_{f_{\psi}}(N(\alpha))| + \ell^{\varepsilon} T^{-3/2} (\ell A)^{-1/2} \sum_{\substack{\alpha \in \mathcal{F}_{D} \\ A/2 \leq N(\alpha) \leq 2A \\ \alpha \notin S_{T^{-1+\delta}} }} |\lambda_{f_{\psi}}(N(\alpha))| \nonumber \\
    &\ll_{D,\varepsilon, \delta} \ell^{\varepsilon} T^{1/2} (\ell A)^{-1/2} \sum_{\alpha\in S_{T^{-1+\delta}}} |\lambda_{f_{\psi}}(N(\alpha))| + \ell^{\varepsilon} T^{-3/2} (\ell A)^{-1/2} \sum_{\substack{\alpha \in \mathcal{F}_{D} \\ A/2 \leq N(\alpha) \leq 2A}} |\lambda_{f_{\psi}}(N(\alpha))|. \nonumber
\end{align}
The result follows from \eqref{BoundInnerSumsSa-1}, \eqref{BoundInnerSumsSa-2}, and \eqref{BoundInnerSumsSa-5}.
\end{proof}

An important observation is that the right hand side of the bound of Proposition \ref{BoundInnerSumsSa} does not depend on $m$, so that we can use well-known bounds for the harmonic sum to estimate $S_a$; for instance, we have 
\begin{equation}
    \sum_{m=1}^{M} \frac{1}{m} \ll \log(M). \label{BoundHarmonicSum}
\end{equation}
Now, from Deligne's proof of the Weil conjectures (see \cite{Deligne-WeilConj}), we have for any $n\geq 1$ that $\lambda_{f_{\psi}}(n) \ll_{\varepsilon} n^{\varepsilon}$. If $\alpha \in \mathcal{F}_{D}$ with $A/2 \leq N(\alpha) \leq 2A$, we have that $$ \lambda_{f_{\psi}}(N(\alpha)) \ll_{\varepsilon} N(\alpha)^{\varepsilon} \ll_{\varepsilon} A^{\varepsilon}.$$ From this bound, Proposition \ref{BoundInnerSumsSa}, and \eqref{BoundHarmonicSum}, we have
\begin{align}
    S_a &\ll_{D,\varepsilon,\delta} \ell^{\varepsilon} T^{1/2} (\ell A)^{-1/2} \sum_{\alpha \in S_{T^{-1+\delta}}} |\lambda_{f_{\psi}}(N(\alpha))| \sum_{m=1}^{(\ell T)^{1+\delta}} \frac{1}{m} \label{MainBound6} \\
    &\qquad + \ell^{\varepsilon} T^{-3/2} (\ell A)^{-1/2} \sum_{\substack{\alpha \in \mathcal{F}_{D} \\ A/2 \leq N(\alpha) \leq 2A}} |\lambda_{f_{\psi}}(N(\alpha))| \sum_{m=1}^{(\ell T)^{1+\delta}} \frac{1}{m} \nonumber \\
    &\ll_{D,\varepsilon,\delta} \ell^{\varepsilon} T^{1/2} (\ell A)^{-1/2} A^{\varepsilon} \log\left((\ell T)^{1+\delta}\right) \left| S_{T^{-1+\delta}} \right| \nonumber \\
    &\qquad + \ell^{\varepsilon} T^{-3/2} (\ell A)^{-1/2} A^{\varepsilon} \log\left((\ell T)^{1+\delta}\right) \sum_{\substack{\alpha \in \mathcal{F}_{D} \\ A/2 \leq N(\alpha) \leq 2A}} 1 \nonumber \\
    &\ll_{D,\varepsilon,\delta} \ell^{\varepsilon} T^{1/2} (\ell A)^{-1/2} (\ell T)^{\varepsilon} \left| S_{T^{-1+\delta}} \right| + \ell^{\varepsilon} T^{-3/2} (\ell A)^{-1/2} A^{\varepsilon} (\ell T)^{\varepsilon} \sum_{\substack{\alpha \in \mathcal{F}_{D} \\ A/2 \leq N(\alpha) \leq 2A}} 1. \nonumber
\end{align}
Now, by a careful analysis of the factorization of rational primes in number fields, one can show for any $\varepsilon>0$ that 
\begin{equation}
    \sum_{N(\mathfrak{a}) = n} 1 \ll_{\varepsilon} n^{\varepsilon}, \label{NumberIdealsGivenNorm}
\end{equation}
where the sum is over all ideals $\mathfrak{a}$ of $\bigo_E$ that have ideal norm equal to $n$. From \eqref{MainBound6} and \eqref{NumberIdealsGivenNorm} it follows that
\begin{align}
    S_a &\ll_{D,\varepsilon,\delta} \ell^{\varepsilon} T^{1/2} (\ell A)^{-1/2} (\ell T)^{\varepsilon} \left| S_{T^{-1+\delta}} \right| + \ell^{\varepsilon} T^{-3/2} (\ell A)^{-1/2} A^{\varepsilon} (\ell T)^{\varepsilon} A^{1+\varepsilon} \label{MainBound7} \\
    &\ll_{D,\varepsilon,\delta} \ell^{2\varepsilon-1/2} T^{\varepsilon + 1/2} A^{-1/2} \left| S_{T^{-1+\delta}} \right| + \ell^{2\varepsilon - 1/2} T^{\varepsilon - 3/2} A^{2\varepsilon + 1/2}. \nonumber
\end{align}
Finally, we estimate the cardinality of the set $S_{R}$ for $0< R \leq 1$.

\subsection{Final estimates}

We will use the Lipschitz principle to estimate $|S_{R}|$. First, since $\mathcal{F}_{D}$ is contained in the right half-plane of $\R^2$, we have that $S_{R} = S_{R'}$ for $R,R'\geq 1/2$, so we only need to bound $|S_{R}|$ for $R\leq 1/2$. Now, consider the set $$ S_{R}^{*} := \left\{ \alpha \in S_{R} : \mathrm{Im}(\alpha)\geq 0 \right\}.$$ Then, we can see that $|S_{R}| \leq 2|S_{R}^{*}|$, so it is enough to estimate $|S_{R}^{*}|$. For this, we will reduce our problem to counting lattice points in a parallelogram. Define the parallelogram $$ \mathcal{P}_{R} := \left\{ av_1 + bv_2 : a,b\in [0,1] \right\},$$ where $$ v_1 := \left( \sqrt{2A}, 0 \right) \text{ and } v_2 := \left( 0, \sin\left(\pi R\right) \sqrt{2A} \right).$$ The important observation is the following.

\begin{lemma}\label{SxiInsideParallelogram}
For $0< R\leq 1/2$, we have $$ S_{R}^{*} \subseteq \mathcal{P}_{R} \cap \bigo_E.$$
\end{lemma}

\begin{proof}
First, note that $S_{R}^{*} \subseteq \mathcal{F}_{D} \subseteq \bigo_E$, so that we only need to show that $S_{R}^{*} \subseteq \mathcal{P}_{R}$. Since $S_{R}^{*}$ is contained in the first quadrant of $\R^2$, we have that $$ S_{R}^{*} = \left\{ \alpha\in \mathcal{F}_{D} : 0\leq \mathrm{Arg}(\alpha) \leq \pi R, \: \frac{A}{2} \leq N(\alpha) \leq 2A \right\}.$$ Suppose $\alpha = x+iy \in S_{R}^{*}$. We then have $A/2 \leq x^2+y^2 \leq 2A$, so that
\begin{equation}
    0\leq x\leq \sqrt{2A}. \label{SxiInsideParallelogram-1}
\end{equation} 
Moreover, since $0\leq \mathrm{Arg}(\alpha) \leq \pi R$ and $0<R\leq 1/2$, then 
\begin{equation}
    0\leq y = \sqrt{x^2+y^2} \sin(\mathrm{Arg}(\alpha)) \leq \sqrt{2A} \sin(\pi R). \label{SxiInsideParallelogram-2}
\end{equation}
The result follows from \eqref{SxiInsideParallelogram-1} and \eqref{SxiInsideParallelogram-2}.
\end{proof}

Hence, to estimate $|S_{R}^{*}|$, we only need to estimate $|\mathcal{P}_{R} \cap \bigo_{E}|$ by Lemma \ref{SxiInsideParallelogram}. Using Proposition \ref{Lipschitz-Dana} with the basis $\{1,i\sqrt{-n}\}$ and the well-known inequality $\sin(x)\leq x$ for $x\geq 0$, we obtain for $0< R\leq 1$ that
\begin{equation}
    |S_{R}| \ll_{D} \left(\sqrt{2A} + R\sqrt{2A}\right) + \left(\sqrt{2A} \cdot R\sqrt{2A}\right) \ll_{D} \sqrt{A} + RA. \label{BoundSetsSxi}
\end{equation}
Combining \eqref{MainBound7} and \eqref{BoundSetsSxi} with $R=T^{-1+\delta}$, and setting $\delta = \varepsilon$, we get
\begin{align}
    S_a &\ll_{D,\varepsilon} \ell^{2\varepsilon-1/2} T^{\varepsilon + 1/2} A^{-1/2} \left( T^{-1+\varepsilon} A + A^{1/2} \right) + \ell^{2\varepsilon -1/2} T^{\varepsilon - 3/2} A^{2\varepsilon + 1/2} \label{MainBound8} \\
    &\ll_{D,\varepsilon} \ell^{2\varepsilon - 1/2} T^{\varepsilon} \left( T^{\varepsilon - 1/2} A^{1/2} + T^{1/2} + T^{-3/2} A^{2\varepsilon + 1/2} \right). \nonumber
\end{align}
However, from the definition of $A$ and $T$, we know that $T\ll \ell$ and $A \ll_{\varepsilon} (T\ell)^{1+\varepsilon}$. Using these estimates in \eqref{MainBound8}, we obtain
\begin{align}
    S_a &\ll_{D,\varepsilon} \ell^{2\varepsilon - 1/2} T^{\varepsilon} \left( T^{\varepsilon - 1/2} (T\ell)^{(1+\varepsilon)/2} + T^{1/2} + T^{-3/2} (T\ell)^{1/2 + 9\varepsilon/2} \right) \label{MainBound9} \\
    &\ll_{D,\varepsilon} \ell^{2\varepsilon - 1/2} T^{\varepsilon} \left( T^{3\varepsilon/2} \ell^{(1+\varepsilon)/2} + T^{1/2} + T^{-1+9\varepsilon/2} \ell^{1/2 + 9\varepsilon/2} \right) \nonumber \\
    &\ll_{D,\varepsilon} \ell^{2\varepsilon - 1/2} \ell^{\varepsilon} \left( \ell^{3\varepsilon/2} \ell^{(1+\varepsilon)/2} + \ell^{1/2} + \ell^{1/2 + 9\varepsilon/2} \right) \ll_{D,\varepsilon} \ell^{8\varepsilon}. \nonumber
\end{align}
Combining \eqref{MainBound3}, \eqref{MainBound4.5}, \eqref{MainBound5}, and \eqref{MainBound9} we obtain the desired result.

\bibliographystyle{alpha}
\bibliography{references}

\end{document}